\PassOptionsToPackage{compatibility=false}{caption}
\documentclass[letterpaper, 10pt, conference, american]{IEEEtran}
\IEEEoverridecommandlockouts
\usepackage[utf8]{inputenc}
\usepackage[american]{babel}
\usepackage{csquotes}
\usepackage[letterpaper,left=48pt,bottom=43pt,right=48pt,top=57pt]{geometry}

\usepackage[
    backend=biber,natbib,
    maxbibnames=99,
    citestyle=ieee,
    bibstyle=ieee,
    sorting=none,
    ]{biblatex}
\usepackage{comment}
\usepackage{url}
\usepackage{graphicx}
\usepackage{subcaption}
\usepackage{multicol}
\usepackage{xspace}
\usepackage{amsmath, amssymb, amsfonts, amsthm}
\usepackage{mathtools}
\usepackage{multirow}
\usepackage{tabularx}
\usepackage[binary-units]{siunitx}

\newcolumntype{C}{>{\centering\arraybackslash}X}
\newcolumntype{L}{>{\raggedright\arraybackslash}X}
\usepackage{siunitx}
\usepackage{listings}
\usepackage{float}

\newcommand{\R}[1][\empty]{\mathbb{R}^{#1}}

\newtheorem{definition}{Definition}
\newtheorem{lemma}{Lemma}
\newtheorem{theorem}{Theorem}

\newtheorem{proposition}{Proposition}[section]
\newtheorem{remarkth}[definition]{Remark}
\newenvironment{remark}{\begin{remarkth}\upshape}{\end{remarkth}}

\newtheorem{corollary}{Corollary}
\newcommand{\bm}[1]{{\boldsymbol{#1}}}

\DeclareMathOperator{\diag}{diag}

\newcommand{\Tr}{\text{Tr}}

\newcommand{\so}{\mathfrak{so}}
\newcommand{\eff}{\text{eff}}

\usepackage{tikz}
\usetikzlibrary{matrix,arrows,calc,positioning,shapes,decorations.pathreplacing}
\tikzstyle{vertex}=[circle,fill=black!20,minimum size=15pt,inner sep=0pt]
\tikzstyle{selected vertex} = [vertex, fill=red!24]
\tikzstyle{edge} = [draw,thick,-]
\tikzstyle{dedge} = [draw,thick,<->]
\tikzstyle{shadowdedge} = [draw, dotted,->]
\tikzstyle{weight} = [font=\small]
\tikzstyle{selected edge} = [draw,line width=5pt,-,red!50]
\tikzstyle{ignored edge} = [draw,line width=5pt,-,black!20]

\usepackage{booktabs}
\begin{document}
\newgeometry{left=48pt,bottom=43pt,right=48pt,top=60pt}

\title{Geometric Control for Load Transportation with Quadrotor UAVs by Elastic Cables}

\author{Jacob R. Goodman$^{1}$, Juan S. Cely$^{1}$, Thomas Beckers$^{2}$ and Leonardo J. Colombo$^{1}$%

%
\thanks{J. Goodman, J. Cely and L. Colombo are with Institute of Mathematical Sciences, ICMAT (CSIC-UAM-UCM-UC3M), Spain. email: {\tt\small \{jacob.goodman, juan.cely, leo.colombo\}@icmat.es}}%

  \thanks{$^{2}$T. Beckers is with Department of Electrical and Systems Engineering, University of Pennsylvania, Philadelphia, USA {\tt\small tbeckers@seas.upenn.edu}}

}

\maketitle
\thispagestyle{empty}
\pagestyle{empty}

\begin{abstract}
Groups of unmanned aerial vehicles (UAVs) are increasingly utilized in transportation task as the combined strength allows to increase the maximum payload. However, the resulting mechanical coupling of the UAVs impose new challenges in terms of the tracking control. Thus, we design a geometric trajectory tracking controller for the cooperative task of four quadrotor UAVs carrying and transporting a rigid body, which is attached to the quadrotors via inflexible elastic cables. The elasticity of the cables together with techniques of singular perturbation allows a reduction in the model to that of a similar model with inelastic cables. In this reduced model, we design a controller such that the position and attitude of the load exponentially converges to a given desired trajectory. We then show that this result leads to an uniformly converging tracking error for the original elastic model under some assumptions. Furthermore, under the presence of unstructured disturbances on the system, we show that the error is ultimately bounded with an arbitrarily small bound. Finally, a simulation illustrates the theoretical results. 
\end{abstract}


\section{Introduction}


The use of aerial robots has become increasingly popular in the last decades due to their superior mobility and versatility in individual and cooperative tasks \parencite{PalCruIRAM12}. For instance, aerial robots can be utilized for monitoring \cite{bayen2006adjoint}, mapping \cite{chung2018survey}, agriculture \cite{ribeiro2021multi}, and delivery \cite{amazonair}, \cite{googlewing}. These vehicles are typically underactuated due to constructional reasons which poses several challenges from the control perspective \cite{rey_dynamics_1999}. Most of the control approaches for individual UAVs are mainly based on feedback linearization \cite{lee_feedback_2009} and backstepping methods \cite{raffo_backsteppingnonlinear_2008} which are also analyzed in terms of stability, e.g. in \cite{frazzoli_trajectory_2000}, \parencite{lee_geometric_2010}.

Multiple aerial robots can be used to transport heavier payloads, thus expanding the capabilities of a single aerial robot \cite{gas}. In the aerial transportation task with multiple quadrotor UAVs considered in this work, a cable establishes a physical connection between the UAV and the cargo. Geometric control of multiple quadrotors with a suspended point-mass load with inelastic cables has been studied in \cite{sreenath_geometric_2013}, and with a rigid body load in \cite{wu_geometric_2014}, \cite{lee_geometric_2018}. In \parencite{god}, the authors model the cables as flexible chains comprised of inflexible links with mass. While these works have considered the cable to be inelastic, we are instead motivated by applications where the elasticity in the cable tethers cannot be ignored without compromising the validity of the estimations \parencite{jo}. In this paper, we study the problem of four quadrotors transporting a rigid body load suspended through elastic cables. Moreover, while the use of elastic cables provide the benefit of reducing impulsive forces on the load, large or rapid oscillations of the load can produce undesired aggressive movements, compromising the load. Therefore, we propose to use elastic cables with high stiffness and damping to guarantee the safety for the rigid load in the transportation task inspired by the fact that most industrial springs are quite stiff and have high damping effects, meaning that larger forces are required to notice the stretching, and oscillations decay quickly. By employing techniques from singular perturbation theory, these assumptions will additionally allow us to reduce the model to that of the corresponding model with inelastic cables, which is differentially flat \cite{sreenath1dynamics}. Such techniques have been applied in the case of a single quadrotor carrying a point-mass load via an inflexible elastic cable in \cite{sreenath_geometric_2017}.

Methods for trajectory tracking and estimation algorithms for underactuated robotic systems evolving on differentiable manifolds such as Lie groups are commonly employed for improving accuracy in simulations. In addition, these methods avoid singularities by working with coordinate-free expressions of controllers  \cite{kobilarov2011discrete}, \cite{kobilarov2014discrete}, \cite{lee2017global}, \cite{izadi2014rigid}. The use of geometric controllers in the UAVs literature has been extensively developed in the last years (see, for instance, \cite{lee_geometric_2010,lee_geometric_2014,lee_geometric_2018,lee_geometric_2013}, and references therein). In this work we advance on the design of geometric controllers for UAVs cooperating to transport a rigid load via inflexible elastic cables. By working directly on the manifold configuration space for the cooperative transportation task between the quadrotors allow us to avoid potential singularities of local parameterizations, e.g. Euler angles, generating agile maneuvers of the payload in a uniform manner. In particular, in this work, a geometric control scheme taking the form of a geometric PD controller (with additional feedback terms) is designed such that both the rigid load's position and attitude track a given desired trajectory exponentially fast. Moreover, we have also considered and theoretically analyzed the situation where unstructured disturbances are applied to the load to evaluate the robustness of the proposed controller. Such a scenario with unstructured disturbances has not been considered before in the literature for point-mass loads nor rigid body loads.

\textbf{Contribution:} The main contributions of this work are: (i) the modeling and subsequent derivation of the corresponding equations of motion for the cooperative task of four quadrotor UAVs transporting a rigid load via inflexible elastic cables. The modeling and dynamics are summarized in Proposition \ref{equationsprop}. (ii) Reduction of these equations of motion to the case of inelastic cables via singular perturbation theory. This can be considered an extension of the results obtained in \cite{sreenath_geometric_2017}, which studies the case of a single quadrotor transporting a point mass load via an inflexible elastic cable. This is developed in Section \ref{sec4} and employed in Section \ref{sec5} where (iii) we provide a geometric control scheme for the exponential tracking of the load position and attitude to some desired trajectories. Lyapunov analysis is used to determine sufficient conditions for exponential tracking. Theorem \ref{Prop: gains} then proves the existence of stabilizing gains—for sufficiently small initial errors in the cables—which satisfy the conditions. Such a proof was not previously seen in the literature, and in principle yields some insight into the relationship between controller gains and exponential stability. Finally, (iv) we handle the case of unstructured bounded disturbances acting on our system, which also had not been seen in the literature. In particular, Theorem \ref{gainsthdist} shows that the same control scheme will yield uniform ultimate bounds in the case of unstructured bounded disturbances acting on the system. Moreover, the ultimate bound can be made arbitrarily small by choosing gains appropriately.

The rest of the paper is structured as follows. In Section \ref{sec3} we model and derive the dynamical system describing the task of carrying and transport a rigid body load between the quadrotors by elastic cables. Section \ref{sec4} reduces the dynamical model introduced in Section \ref{sec3} by employing singular perturbation theory techniques. The main results of this work showing the convergence of the tracking error for the reduced and the actual system of quadrotors are given in Section \ref{sec5} and \ref{sec_dist}. In Section \ref{sec_dist}, we introduce unstructured bounded disturbances to the reduced model and show that the same control scheme can be applied to achieve uniform ultimate boundedness. Finally, numerical simulations visualize the theoretical results. 

\section{Modelization and Control Equations}\label{sec3}
We begin by modelling and deriving the control system describing the transportation task between the quadrotors. This is done by constructing the total kinetic and potential energies of the system, in addition to the virtual work done by non-conservative forces and, subsequently, by using tools from variational calculus on manifolds \parencite{Leebook}, \parencite{HSS}. 

\begin{figure}[htb]
\centering
\includegraphics[width=8cm]{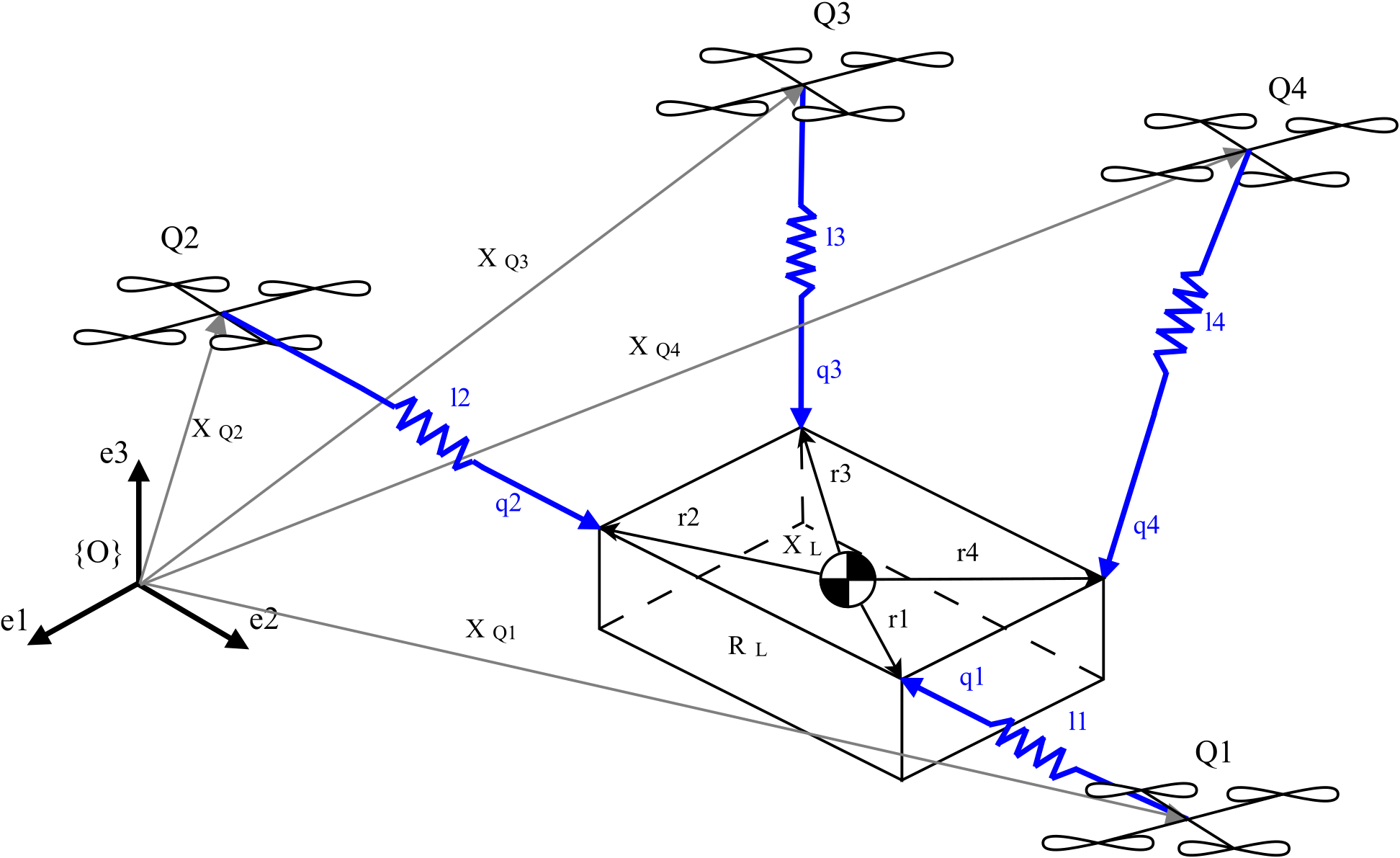}
\caption{Proposed model for the transportation task with quadrotor UAVs.} 
\label{fig:drone_sys}
\end{figure}

Consider four identical quadrotor UAVs transporting a rigid body of total mass $m_L\in\mathbb{R}_{>0}$ and positive-definite inertia matrix $J_L\in\R[3\times 3]$. The load is considered rigid and of uniform mass density. It is connected to the center of mass of each quadrotor via a massless inflexible elastic cable of rest length $L$. A graphical description of the proposed system is visible in Figure \ref{fig:drone_sys}.

The configuration space of the mechanical system describing the cooperative task between quadrotors is given by $$Q = \underbrace{(SO(3) \times \R[3])}_\text{Rigid load} \times \underbrace{(S^2 \times \R)^4}_\text{Cables}\times\underbrace{(SO(3))^4}_\text{Quadrotor attitudes},$$ where $SO(3)$ denotes the special orthogonal group of $3\times 3$ rotation matrices and $S^2$ denotes the $2$-sphere. 
The basic notation and methodology along this section is fairly
standard within the geometric control and classical mechanics literature, and we have attempted to use traditional symbols and definitions wherever feasible. Table~\ref{tbl:nomenclature} provides the symbols and geometric spaces that are used frequently throughout the paper. 

The system has $30$ degrees of freedom—$6$ degrees corresponding to the load, $3$ degrees for each cable, and $3$ degrees for each quadrotor. Observe that the positions of the quadrotors do not appear in the configuration space, as they are uniquely defined in terms of the other state variables due to the constraints $x_{Q_j} = x_L + R_Lr_j - l_j q_j$ for $j \in \{1,2,3,4\}$. Meanwhile, we have $16$ inputs to the system in the form of four thrust controls $f_j \in \R$, corresponding to the total lift force exerted on the quadrotors by the spinning propellers, and four moment controls $M_j \in \R[3]$, which are related to the torques induced on the quadrotors by the rotating propellers—for each $j \in \{1,2,3,4\}$. Thus the complete systems has $14$ degrees of under-actuation, with the quadrotor positions and attitudes being directly actuated. Note also that upon fixing a body frame to each quadrotor such that the vector $\bar{e}_3 = [0, 0, 1]^T$ points in the direction of the applied thrust, we may alternatively express the thrust controls $f_j$ via the vectors $u_j=f_jR_je_3\in\mathbb{R}^3$ in the inertial frame for $j \in \{1,2,3,4\}$. We will utilize this representation frequently throughout the paper. Alternatively, one may choose to control the total thrust of \textit{each} propeller individually. However we opt for the former approach both because it is more pervasive in the literature and because it leads nicely to the separation of the quadrotor's attitude dynamics from the rest of the system's dynamics.  The thrust generated by the $i$-th propeller along the $e_3$-axis can be determined by the total thrust and the moment controller \cite{lee_geometric_2010}.

\begin{table}[!htb]
\centering
\begin{tabularx}{\columnwidth}{ccX}
\toprule
\textbf{Symb.} & \textbf{Space} & \textbf{Description} \\
\midrule
$x_L$             & $\R[3]$        & Position of the center of mass of the load in inertial frame. \\
$v_L$             & $\R[3]$        & Translational velocity of center of mass of the load in the inertial frame. \\
$R_L$             & $SO(3)$          & Attitude of the load in body frame. \\
$\omega_L$        & $TSO(3)$         & Angular velocity of the load in inertial frame. \\
$m_L$             & $\R$           & Mass of the load. \\
$J_L$             & $\hbox{Sym}_{\succ 0}( \R)$ & Moment of inertia of the load. \\
\midrule
$x_{Q_j}$         & $\R[3]$        & Position of quadrotor $j$ in inertial frame. \\
$R_j$             & $\text{SO}(3)$ & Attitude of quadrotor $j$. \\
$\Omega_j$        & $\so(3)$       & Angular velocity of quadrotor $j$ in body frame. \\
$m_{Q}$         & $\R$           & Mass of quadrotors. \\
$J_{Q}$         & $\hbox{Sym}_{\succ 0}(\R)$    & Moment of inertia of quadrotors. \\
$u_j$             & $\R[3]$        & Net thrust applied vertically in the body frame of quadrotor $j$. \\
$M_j$             & $\R[3]$        & Moment vector in body frame of quadrotor $j$. \\
\midrule
$L$             & $\R$           & Rest length for elastic cables. \\
$r_j$             & $\R[3]$        & Unit vector from the center of mass of the load to quadrotor $j$. \\
$q_j$             & $S^2$          & Position vector of cable suspended from quadrotor $j$. \\
$l_j$             & $\R$           & Length of elastic cable attached to quadrotor $j$. \\
$\omega_j$        & $TS^2$         & Angular velocity of cable $j$ in the inertial frame. \\
$u^{\perp_j}$     & $\R[3]$        & The component of $u$ that is perpendicular to $q_j$. \\
$u^{\parallel_j}$ & $\R[3]$        & The component of $u$ that is parallel to $q_j$. \\
\bottomrule
\end{tabularx}
\caption{Nomenclature}
\label{tbl:nomenclature}
\end{table}

The translational kinetic energy of each quadrotor can be described by $\frac12 m_Q ||\dot{x}_{Q_j}||^2$, where $m_Q$ denotes the mass of the quadrotor. As the quadrotors and load are rigid bodies, we further have rotational kinetic energy components in the total kinetic energy. Fixing a body frame to each quadrotor and denoting the angular velocity in this body frame by $\Omega_j\in\R[3]$, the angular kinetic energy is given by $\frac{1}{2}\Omega_j^T J_Q \Omega_j$, where $J_Q$ is a symmetric positive-definite inertia tensor. The angular velocity $\Omega_j$ is defined by the kinematic equation $\dot{R}_j = R_j \hat{\Omega}_j$, where $\hat{\cdot}: \R[3] \to \so(3)$ is the \textit{hat isomorphism} which maps vectors on $\R[3]$ to the set $\mathfrak{so}(3)$ of $(3\times 3)$ skew-symmetric matrices, that is, 
$$\displaystyle{\Omega = \begin{bmatrix} \Omega_1 \\ \Omega_2 \\ \Omega_3 \end{bmatrix} \mapsto \begin{bmatrix} 0 & -\Omega_3 & \Omega_2 \\ \Omega_3 & 0 & -\Omega_1 \\ -\Omega_2 & \Omega_1 & 0 \end{bmatrix} := \hat{\Omega}}.$$ 

We follow the same procedure for the load to obtain the angular kinetic energy as $\frac12 \Omega_L^T J_L \Omega_L$, where $J_L$ is the inertia tensor of the load (again a symmetric positive-definite matrix) and $\Omega_L$ is the angular velocity defined by $\dot{R}_L = R_L \hat{\Omega}_L$.

Elements in the tangent space $T_{R}SO(3)$ at $R\in SO(3)$ are identified with elements in $SO(3)\times\mathfrak{so}(3)$ by a left-trivialization. That is, for $R\in SO(3)$, the map $(R,\dot{R})\in T_{R}SO(3)\mapsto (R,R^{-1}\dot{R})=:(R,\hat{\Omega})\in SO(3)\times\mathfrak{so}(3)$ is a diffeomorphism (see \parencite{HSS} for details). Therefore, after a left trivialization of $TSO(3)$, the tangent bundle of $Q$, describing the state space of the system, can be identified as
$TQ \cong  (SO(3)\times\so(3))^4\times(S^2\times TS^{2}\times\R\times\R)^4\times(SO(3)\times\so(3)\times \R[3]\times \R[3])$, where we have used that the tangent bundle of a finite dimensional vector space $V$, i.e, $TV$, is isomorphic to $V\times V$.

Finally, the total kinetic energy of the system, $K: TQ \to \R$,  is given by summing the respective translational and angular kinetic energies of the quadrotors and load, that is \begin{equation*}
\begin{split}
K & = \underbrace{\frac{1}{2}m_L ||\dot{x}_L||^2 + \sum_{j=1}^{4} \frac{1}{2}m_Q ||\dot{x}_{Q_j}||^2}_\text{Translational K.E.} \\
  & \hspace{2.5cm}+ \underbrace{\frac12 \Omega_L^T J_L \Omega_L + \sum_{j=1}^{4} \frac{1}{2}\Omega_j^T J_Q \Omega_j. }_\text{Angular K.E.}
\end{split}
\end{equation*}

Moreover, the total potential energy of the system, $U: Q \to \R$, is given by $$U = \sum_{j=1}^{4}(\underbrace{m_Q ge_3^T x_{Q_j} + m_L g e_3^T x_L}_\text{Gravitational P.E.} +  \underbrace{ \frac{1}{2}k(L - l_j)^2}_\text{Elastic P.E.}),$$which corresponds to the gravitational potential energies of the quadrotors and the load, as well as the elastic potential of the cables. Observe that here we have fixed an inertial frame such that $e_3$ is oriented opposite to the direction of the gravitational acceleration. As usual, the Lagrangian of the system $\mathbf{L}: TQ \to \R$ is defined by $\mathbf{L}:= K - U$.

Note that the control inputs take the form of non-conservative external forces, so that we must use Lagrange d'Alembert Variational Principle (see \parencite{abloch} for instance) -- with controls as virtual forces in our system -- to obtain our system dynamics from the Lagrangian $\mathbf{L}$. We further wish to add a non-conservative force corresponding to a damping in the elastic cables. That is, a velocity dependent force that serves to reduce the amplitude of oscillations in the elastic cables. In particular, we will opt to make this force \textit{proportional} to the velocity—as is standard for damped harmonic oscillators—with constant of proportionality $c > 0$. 

Denote by $C^{\infty}([0,T],Q,q_0,q_T)$ the space of smooth function from $[0,T]$ to $Q$ with fixed endpoints points, denoted by $q_0$ and $q_T$, respectively. Consider the action functional $\mathcal{A}:C^{\infty}([0,T],Q,q_0,q_T)\to\R$ given by 
\begin{equation*}
\begin{split}
    \mathcal{A}(c(t))&=\int_{0}^{T} \mathbf{L}(c(t),\dot{c}(t))\, dt \\
    &+ \sum_{j=1}^{4} \int_{0}^{T} \left(||f_jR_je_3||_{\R[3]}^{2} + ||\hat{M}_j||^2_{\mathfrak{so}(3)} - c \dot{l}_j\right)\,\, dt,
\end{split}
\end{equation*}  
where $||\hat{M}_j||_{\mathfrak{so}(3)}:=\langle\hat{M}_j,\hat{M}_j\rangle^{1/2}=\sqrt{\hbox{Tr}(\hat{M}_j^{T}\hat{M}_j)}$, with  $c(t):=(R_L(t),x_L(t),q_j(t),l_j(t),R_j(t)))\in C^{\infty}([0,T],Q,q_0,q_T)$. It is well know (see \cite{MR}, \cite{Leebook}, \parencite{HSS} for instance) that critical points of the action functional $\mathcal{A}$ corresponds with the forced Euler-Lagrange equations for the Lagrangian $\mathbf{L}$ and external forces described before by the control inputs and damping force.

In order to use Lagrange-d'Alembert principle, we must describe the variations of our state variables. These variations must be tangent vectors in the tangent spaces of the submanifolds of the configuration space in which the state variables live. In addition, they must vanish at the end points, because tangent vectors on the tangent bundle of $C^{\infty}([0,T],Q,q_0,q_T)$ must satisfy such a condition (see for instance \cite{Leebook}, \parencite{HSS}, \parencite{MR}). 

In particular, we choose $\delta x_L \in \R[3]$ and $\delta l_j \in \R$  arbitrary,  $\delta q_j = \frac{d}{d\epsilon}\mid_{\epsilon=0}\exp(\epsilon \hat{\xi}_j)q_j=\xi_j \times q_j \in T_{q_j} S^2$, satisfying $\xi_j\cdot q_j=0$ for arbitrary vectors $\xi_j \in \R[3]$ and $j = 1,...,4$. By additionally defining the curves on the Lie algebra $\mathfrak{so}(3)$ given by $\hat{\eta}_j = R_j^T \delta R_j \in \so(3)$, it can be shown that (see for instance \parencite{MR} Chapter 13) $\widehat{\delta \Omega}_j = \widehat{\hat{\Omega}_j \eta_j} + \dot{\hat{\eta}}_j$ with $\hat{\eta}_j$ satisfying $ \hat{\eta}_j(0)= \hat{\eta}_j(T)=0$ (since $\delta R_j(0)=\delta R_j(T)=0)$ for $j = 1,...,4$. Moreover, we have the following relations
\begin{align*}
\delta x_{Q_j} &= \delta x_L + \delta R_L r_j - \delta l_j q_j - l_j \delta q_j, \\
\delta \dot{x}_{Q_j} &= \delta \dot{x}_L + \delta \dot{R}_L r_j - \delta \dot{l}_i q_i - \dot{l}_i \delta q_i - \delta l_i \dot{q}_i - l_i \delta \dot{q}_i.
\end{align*}

The controlled dynamics of the system is described as follow by finding critical points of $\mathcal{A}$ by employing variational calculus on differentiable manifolds.

\begin{proposition}\label{equationsprop}
Critical points of the action functional $\mathcal{A}$ for variations with fixed endpoints corresponds with solutions of the controlled Euler-Lagrange equations
\begin{align}
&\dot{x}_L = v_L,\quad \dot{R}_L = R_L \hat{\Omega}_L, \label{xkin} \\ 
&m_{\eff} ( \dot{v}_L + ge_3) = \sum_{j=1}^4 u_j + m_Q \ddot{\zeta}_j - m_Q R_L(\hat{\Omega}_L^2 + \dot{\hat{\Omega}}_L)r_j,\label{xdyn}\\ 
 &J_{\eff} \dot{\Omega}_L + \hat{\Omega}_L J_{\eff} \Omega_L = \sum_{j=1}^4 m_Q \hat{r}_j R_L^T (-g e_3 - \dot{v}_L \label{RLdyn}\\& \hspace{5.9cm}+ \ddot{\zeta}_j + \frac1{m_Q}u_j), \nonumber \\ 
&m_Q q_j^T \ddot{\zeta}_j = m_Q q_j^T ( \dot{v}_L + R_L(\hat{\Omega}_L^2 + \dot{\hat{\Omega}}_L)r_j + ge_3 - \frac1{m_Q}u_j) \nonumber\\& - c \dot{l}_j  + k(L - l_j), \label{ldyn}\\ 
&q_j \times \ddot{\zeta}_j = q_j \times ( \dot{v}_L + R_L(\hat{\Omega}_L^2 + \dot{\hat{\Omega}}_L)r_j + ge_3 -\frac1{m_Q}u_j)\label{qjdyn} \\ 
&J_Q \dot{\Omega}_j = J_Q \Omega_j \times \Omega_j + M_j,\,\,\dot{R}_j = R_j \hat{\Omega}_j,\,\, \ j = 1,\ldots,4 \label{Rdyn}
\end{align} where $m_{\eff} := 4m_Q + m_L, \ J_{\eff} := J_L - \sum_{j=1}^4 m_Q \hat{r}_j^2$, and $\zeta_{j} := l_j q_j$.
\end{proposition}
\textit{Proof:} The proof follows similar arguments than proof of Proposition $1$ in \cite{jacobquads} and in \cite{lee_geometric_2018} by expanding the previous variations in the action functional and applying the Fundamental Lemma of Calculus of Variations \parencite{GF}. We include the detailed proof in the Appendix.\hfill $\square$
\begin{remark}
Equations \eqref{xkin} describe the kinematics of the load's position and attitude, respectively.  Similarly, equations \eqref{xdyn}-\eqref{RLdyn} describe the dynamics of the load's position and attitude, respectively, and equations \eqref{Rdyn}, corresponds with the dynamics and kinematics of the attitudes of each quadrotor, respectively. 

Equations \eqref{ldyn}, indexed for $j = 1,\ldots,4$, describe the dynamics of the lengths of the elastic cables. This can be understood by observing that the projection of $\ddot{\zeta}_j$ onto $q_j$ preserves the acceleration of the length (which is inherently oriented along the cable), while removing the acceleration of the attitude from consideration with the identity $q_j^T \ddot{q}_j = -||\dot{q}_j||^2$. Conversely, equation \eqref{qjdyn}, indexed for $j = 1,\ldots,4$, describe the dynamics of the attitudes of the elastic cables, as the cross-product with $q_j$ preserves the acceleration of the cable attitude while annihilating the acceleration of the cable length. 
\end{remark}


\begin{remark}
Note that the proposed dynamics given by Proposition \ref{equationsprop} can be easily extended to an arbitrary number $n$ of quadrotors and elastic cables by making the range of sum in the kinetic and potential energies evolve from $1$ to $n$ instead of only from $1$ to $4$. We choose $n=4$ only for illustrative purposes on the transportation task, but all the results in this paper follows exactly the same procedure for the case of $n$ by appropriately change the range of the sums. Note that the analysis conducted in \cite{lee_geometric_2018} does not incorporate the elasticity of the cables which is the main difference with respect to the model proposed for our system. In addition one can consider different kind of objects. The difficulty here is the inertia mass of the load. In order to apply the results of our work with the proposed control strategy the inertia mass of the load $J_L$ must be a non-singular matrix. For instance, in the case of transport a rigid bar instead of a load, $J_L$ is singular and the control design must be conducted in a different way as it was shown in our previous paper \cite{jacobquads}.

\end{remark}

\section{Reduced Model}\label{sec4}
While the use of elastic cables provide the benefit of reducing impulsive forces on the load, large or rapid oscillations of the load can produce undesired aggressive movements, compromising the load. This can be combated by utilizing elastic cables with high stiffness and damping to guarantee the safety for the rigid load in the transportation task. Such a condition will commonly be fulfilled by the cables used in applications. 

A benefit of this assumption is that we will be able to reduce the degrees of freedom in the original model via techniques from singular perturbation theory \cite{SPT}. In fact, we will see that the reduced model is \textit{differentially flat}. In other words, the states and inputs of the reduced model can be written as algebraic functions of $16$ \textit{flat outputs} and their derivatives—which dramatically reduces the difficulty involved in generating dynamically feasible trajectories for under-actuated systems. On the other hand, the original elastic model is not differentially flat, which further justifies our desire to reduce the model.

In particular, we will consider the case that $\displaystyle{k = \frac{\bar{k}}{\epsilon^2}}$ and $\displaystyle{c = \frac{\bar{c}}{\epsilon}}$ with $\bar{k}, \bar{c} > 0$ and $\epsilon > 0$ sufficiently small, and we will show that the dynamics approach that of the same model with inelastic cables (that is, with $l \equiv L$) as $\epsilon \to 0$. We further consider a change of variables of the form $l_j = \epsilon^2 y_j + L$ and $\dot{l}_j = \epsilon z_j$, which is motivated by observing that $k( L - l_j) = \bar{k}y_j$ and $c \dot{l}_j = \bar{c} z_j$. From this, we can see that $\zeta_j = (\epsilon^2 y_j + L)q_j$. Therefore, $\ddot{\zeta}_j = L\ddot{q}_j + \epsilon (\dot{z}_j q_j + z_j \dot{q}_j) + \epsilon^2 y_j \ddot{q}_j$. Making these substitutions into the dynamics described in Proposition \ref{equationsprop}, in addition to employing the equation $\dot{q}_j = \omega_j \times q_j$, we obtain
\begin{align*}
&\dot{x}_L = v_L,\,\,\dot{R}_L= R_L \hat{\Omega},\,\,\epsilon \dot{y}_j = z_j,\,\, \dot{q}_j = \omega_j \times q_j, \\ 
&m_{\eff} ( \dot{v}_L + ge_3) = \sum_{j=1}^4 \Big(u_j - m_Q R_L(\hat{\Omega}_L^2 + \dot{\hat{\Omega}}_L)r_j \\
&\hspace{2.5cm}+m_Q (L\ddot{q}_j + \epsilon (\dot{z}_j q_j + z_j \dot{q}_j) + \epsilon^2 y_j \ddot{q}_j) \Big),\\ 
&J_{\eff} \dot{\Omega}_L + \hat{\Omega}_L J_{\eff} \Omega_L = \sum_{j=1}^4 m_Q\hat{r}_j R_L^T \Big(-g e_3 - \dot{v}_L + L\ddot{q}_j \\ 
&\hspace{2.8cm} +\epsilon (\dot{z}_j q_j + z_j \dot{q}_j) + \epsilon^2 y_j \ddot{q}_j + \frac1{m_Q}u_j\Big),\,\,%
\end{align*} 
\begin{align*}
&m_Q L \epsilon \dot{z}_j =q_j^T \Big( \bar{c}z_j q_j + \bar{k} y_j q_j + m_Q u_j + m_Q L(1  + \epsilon^2 y_j)\ddot{q}_j \\
&\hspace{2cm}- m_Q ( \dot{v}_L - m_Q R_L(\hat{\Omega}_L^2 + \dot{\hat{\Omega}}_L)r_j+ ge_3) \Big),  \\ 
&\dot{\omega}_j = \frac1{L} q_j \times  \Big( \dot{v}_L - R_L(\hat{\Omega}_L^2 + \dot{\hat{\Omega}}_L)r_j + ge_3 \\
&\hspace{2cm}- \frac1{m_Q} u_j - \epsilon ( \dot{z}_j q_j + z_j \omega_j) - \epsilon^2 y_j \ddot{q}_j \Big), \\
&J_Q \dot{\Omega}_j = J_Q \Omega_j \times \Omega_j + M_j,\,\, \dot{R}_j = R_j \hat{\Omega}_j,\,\, j = 1,...,4.
\end{align*} 

The previous system can be written as \begin{equation}\label{spmodel}
\dot{x} = f(t, x, z; \epsilon), \,\,\,\, \epsilon \dot{z} = g(t, x, z; \epsilon),\end{equation}where $f$ and $g$ are smooth functions, $x$ is the vector representing $(x_L, v_L, R_L, \Omega_L, q_j, \omega_j, R_j, \Omega_j)$, and $z$ is the vector representing $(y_j, z_j)$, for $j = 1,\ldots,4$. 

The above dynamical system is known as singular perturbation model \cite{SPT}, with the first equation describing the slow dynamics and the second describing the fast dynamics. Evaluating at $\epsilon = 0$, the fast dynamics provide us with algebraic equations that can be solved to obtain $z = h(t, x)$. In particular, \begin{align}
z_j &= 0, \label{eqh}\\
y_j &= -\frac{m_Q}{\bar{k}} q_j^T \left[ u_j + L\ddot{q}_j - \dot{v}_L + R_L(\hat{\Omega}_L^2 + \dot{\hat{\Omega}}_L)r_j- ge_3 \right].
\end{align}

Substituting these equations back into the slow dynamics, we obtain the reduced (slow) model of the control system describing the cooperative task, given by $\dot{x} = f(t, x, h(t, x), 0)$. That is,
\begin{align}
&\dot{x}_L = v_L,\,\,\dot{R}_L = R_L \hat{\Omega},\,\,\dot{q}_j = \omega_j \times q_j, \label{xred} \\ 
&m_{\eff} ( \dot{v}_L + ge_3) = \nonumber\\&\hspace{1cm}\sum_{j=1}^4 \Big(u_j - m_Q R_L(\hat{\Omega}_L^2 + \dot{\hat{\Omega}}_L)r_j + m_Q L\ddot{q}_j \Big), \label{vred}\\ 
&J_{\eff} \dot{\Omega}_L + \hat{\Omega}_L J_{\eff} \Omega_L =\nonumber\\ &\hspace{1cm}\sum_{j=1}^4 m_Q \hat{r}_j R_L^T \Big(-g e_3 - \dot{v}_L + L\ddot{q}_j + \frac1{m_Q}u_j\Big),\label{OmegaLred} \\ 
&\dot{\omega}_j = \frac1{L} \hat{q}_j\Big( \dot{v}_L - R_L(\hat{\Omega}_L^2 + \dot{\hat{\Omega}}_L)r_j + ge_3 - \frac1{m_Q} u_j \Big), \label{wjred} \\
&J_Q \dot{\Omega}_j = J_Q \Omega_j \times \Omega_j + M_j,\,\, \dot{R}_j = R_j \hat{\Omega}_j,\,\, j = 1,...,4 \label{Rdynred}.
\end{align} 
\begin{remark}
Observe that the reduced model \eqref{xred} - \eqref{Rdynred} preserves the original 16 inputs of the system, but the configuration space has lost 4 degrees of freedom (namely the cable lengths). Hence, the reduced model has $10$ degrees of under-actuation. In fact, it can be seen that the reduced model is equivalent to the original model with inelastic cables \cite{lee_geometric_2018} (that is, where $l \equiv L$). This system was shown to be differentially flat in \cite{sreenath1dynamics}. From Proposition $1$ in \cite{lee_geometric_2018}, we know that achieving exponentially stable tracking of the reduced model on some set of initial conditions will guarantee exponentially stable tracking in some subset of those initial conditions—whose relative size depends on $\epsilon$. Hence, we may work within this reduced model to design geometric controllers towards the end of tracking the position and attitude of the load.\end{remark}


\section{Control Design for Position and Attitude Trajectory Tracking of the Rigid Load}\label{sec5}

\begin{figure}[htb]
\centering
\includegraphics[width=5.5cm]{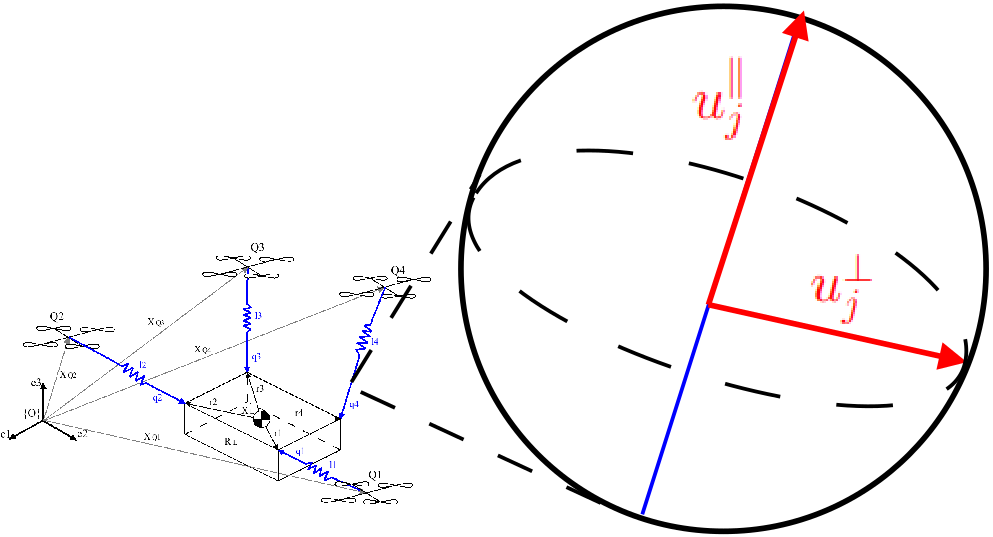}
\caption{Graphical description for the decomposition of the control input viewed on each cable.} 
\label{fig:control_inputs}
\end{figure}

In the following we discuss the control design for the reduced model \eqref{xred}-\eqref{Rdynred}. That is, we provide a set of controllers $u_j\in\mathbb{R}^{3}$ such that the position and attitude of the load reach a desired position $\tilde{x}_L\in\mathbb{R}^{3}$ and attitude $\tilde{R}_L\in SO(3)$ exponentially fast. The strategy is to decompose $u_j$ into components which are parallel and perpendicular to the cable attitudes $q_j$ via $u_j = u_j^{\parallel} + u_j^{\perp}$, where $u_j^{\parallel} = (q_j^T u_j)q_j$ and $u_j^{\perp} = (I - q_j q_j^T)u_j$. This is motivated by the fact that only the $u_j^{\parallel}$ components influence the dynamics of the load, while only the $u_j^{\perp}$ components influence the cables dynamics. Its graphical description is shown in Figure \ref{fig:control_inputs}. 

In particular, we show that the reduced model is equivalent to the dynamical system discussed in \cite{lee_geometric_2018} with $n = 4$ quadrotors and $l_j = L$ for $j = 1,...,4$. As such, we may implore the control scheme designed there. After this, we will introduce configuration error functions for each state variable, from which we may derive our error dynamics.

Notice that equations \eqref{Rdynred} for quadrotor attitude are independent from the rest of the dynamics, and the moment controllers $M_j$ appear exclusively within them. Moreover, these equations are just those that appear in \cite{lee_geometric_2010} and \cite{lee_geometric_2018}, for which $M_j$ was designed to attain almost-global exponential stability. We use the same controller for the attitude of the quadrotors, and disregard the equations for the remainder of this paper.

\subsection{Error dynamics and control design}
We begin by further simplifying the dynamical system \eqref{xred}-\eqref{wjred}. In particular, we find an equation for $\ddot{q}_j$ that we will substitute into equations \eqref{vred} and \eqref{OmegaLred}. By differentiating $\dot{q}_j = \omega_j \times q_j$ and expanding it with the vector triple product identity, it can be shown that $\ddot{q}_j = \dot{\omega}_j \times q_j - ||\omega_j ||^2 q_j$. Now we may substitute \eqref{wjred} in for $\dot{\omega}_j$ to find that

\begin{align*}
L \ddot{q}_j &= -q_j \times (L\dot{\omega}_j) - L||\omega_j||^2 q_j \\
&= (I - q_j q_j^T)\Big[\dot{v}_L - R_L(\hat{\Omega}_L^2 + \dot{\hat{\Omega}}_L)r_j + ge_3 - \frac1{m_Q} u_j \Big] \\
&\quad -L\left|\left|\omega_j\right|\right|^2 q_j.
\end{align*} 
Substituting this equation for $m_Q L \ddot{q}_j$ into \eqref{vred} and making use of the fact that $m_{\eff} = 4m_Q + m_L$, we obtain
\begin{align}
M_L (\dot{v}_L + ge_3) = &\sum_{j=1}^4 \left[m_Qq_j q_j^T R_L (\hat{\Omega}_L^2 + \dot{\hat{\Omega}}_L)r_j\right.\nonumber\\ &+\left. u_j^{\parallel} - m_Q L \left|\left|\omega_j\right|\right|^2 q_j \right],\label{xLdyn_Lee}
\end{align}where $M_L = m_L I + \sum_{j=1}^4 m_Q q_j q_j^T$. Repeating this procedure with \eqref{OmegaLred}, and making use of the fact that $J_{\eff} = J_L - \sum_{j=1}^4 m_Q r_j$ yields
\begin{align}
&J_L \dot{\Omega}_L + \hat{\Omega}_L J_L \Omega_L= \sum_{j=1}^4 m_Q \hat{r}_j R_L^T\Big[ q_j q_j^T R_L(\hat{r}_j \dot{\Omega}_L + \hat{\Omega}_L^2 r_j) \nonumber\\
&\hspace{1cm}- q_j q_j^T(\dot{v}_L + ge_3) - L||\omega_j||^2q_j + \frac{1}{m_Q}u_j^{\parallel}\Big]\label{OmegaLdyn_Lee}.
\end{align}

With equation \eqref{xLdyn_Lee} in place of \eqref{xred}, and \eqref{OmegaLdyn_Lee} in place of \eqref{OmegaLred}, the reduced model is equivalent to the dynamical system discussed in \cite{lee_geometric_2018} with $n = 4$ and $l_j = L$ for $j = 1,...,4.$. In particular, if we denote our desired trajectories by $\tilde{x}_L, \tilde{v}_L, \tilde{R}_L, \tilde{\Omega}_L, \tilde{q}_j, \tilde{\omega}_j$ and define tracking errors as: \begin{align*}e_{x_L} &= x_L - \tilde{x}_L,\\ e_{v_L} &= v_L - \tilde{v}_L,\\ 
    e_{R_L} &= \frac12 \left(\tilde{R}_L^T R_L - R_L^T \tilde{R}_L \right)^{\vee},\\ e_{\Omega_L} &= \Omega_L - R_L^T \tilde{R}_L \tilde{\Omega}_L,\\
    e_{q_j} &= \tilde{q}_j \times q_j,\\ 
    e_{\omega_j} &= \omega_j + \hat{q}_j^2 \tilde{\omega}_j,\end{align*} then for some choice of gains $k_{x_L}, k_{v_L}, k_{R_L}, k_{\Omega_L}, k_q, k_{\omega} \in \R$, the controllers take form:
\begin{align}
    u_j^{\parallel} &= \mu_j + m_Q L ||\omega_j||^2 q_j + m_Q q_j q_j^T a_j, \label{controllers}\\
    u_j^{\perp} &= m_Q L \hat{q}_j \left(-k_q e_{q_j} - k_{\omega} e_{\omega_j} - (q_j^T \tilde{\omega}_j)\dot{q}_j - \hat{q}_j^2 \dot{\hat{\omega}}\right) \label{controllers2}\\&- m_Q \hat{q}_j^2 a_j,\nonumber
\end{align}where $a_j = \dot{v}_L + ge_3 + R_L \hat{\Omega}^2_L r_j - R_L \hat{r}_j \dot{\Omega}_L$, and $\mu_j$ is an additional controller which satisfies:
\begin{align}
    \mu_j &= q_j q_j^T \tilde{\mu}_j\label{cont1} \\
    \sum_{j=1}^4 \tilde{\mu}_j &= m_L\left(-k_{x_L} e_{x_L} - k_{v_L} e_{v_L} + \dot{\tilde{v}}_L + ge_3 \right) \label{cont2}\\
    \sum_{j=1}^4 \hat{r}_j R_L^T \tilde{\mu}_j &= -k_{R_L} e_{R_L} - k_{\Omega_L} e_{\Omega_L} +\label{cont3} \\
    & \quad \left(R_L^T \tilde{R}_L \tilde{\Omega}_L \right)^{\wedge} J_LR_L^T \tilde{R}_L \tilde{\Omega}_L + J_L R_L^T \tilde{R}_L \dot{\tilde{\Omega}}_L.\nonumber
\end{align}
We further define the desired cable attitudes $\tilde{q}_j = -\frac{\tilde{\mu}_j}{||\tilde{\mu}_j||}$, so that we have $\tilde{\mu}_j \to \mu_j$ as $\tilde{q}_j \to q_j$. 

Notice that equations \eqref{cont2} and \eqref{cont3} can also be written in the form
$$\mathcal{P} \text{diag}(R_L^T,..., R_L^T) \begin{bmatrix} \tilde{\mu}_1 \\ \vdots \\ \tilde{\mu}_4 \end{bmatrix} = \begin{bmatrix} R_L^T \tilde{F} \\ \tilde{M} \end{bmatrix}$$
where
\begin{align*}
\mathcal{P} &= \begin{bmatrix} I_{3\times 3} & \cdots & I_{3 \times 3} \\ \hat{r}_1 & \cdots & \hat{r}_4 \end{bmatrix} \in \R[6\times 12],\\ \tilde{F} &= m_L\left(-k_{x_L} e_{x_L} - k_{v_L} e_{v_L} + \dot{\tilde{v}}_L + ge_3 \right),\\ 
\tilde{M}&=-k_{R_L} e_{R_L} - k_{\Omega_L} e_{\Omega_L} + \left(R_L^T \tilde{R}_L \tilde{\Omega}_L \right)^{\wedge} J_L R_L^t \tilde{R}_L \tilde{\Omega}_L\\ &+ J_L R_L^T \tilde{R}_L \dot{\tilde{\Omega}}_L.
\end{align*}
Then, if we assume that $\text{rank}(\mathcal{P}) \ge 6$ (which depends on the physical connection points of the cables to the load), there is guaranteed to exist a solution $\begin{bmatrix}\tilde{\mu}_1 & \cdots & \tilde{\mu}_4\end{bmatrix}^T$. However, in general there will exist multiple solutions, so we choose the solution with minimal (Euclidean) norm. This is given by
\begin{align}
     \begin{bmatrix} \tilde{\mu}_1 \\ \vdots \\ \tilde{\mu}_4 \end{bmatrix} = \text{diag}(R_L, ..., R_L) \mathcal{P}^T (\mathcal{P} \mathcal{P}^T)^{-1} \begin{bmatrix} R_L^T \tilde{F} \\ \tilde{M} \end{bmatrix}.
\end{align} With these controllers, the error dynamics take the following form
\begin{align}
    &\dot{e}_{v_L} = -k_{v_L}e_{v_L} - k_{x_L}e_{x_L} - \frac1{m_Q} \sum_{j=1}^4 (\tilde{q}_j^T \tilde{\mu}_j)\hat{q}_j e_{q_j},\label{xerror} \\
    &J_L \dot{e}_{\Omega_L} = -k_{\Omega_L} e_{\Omega_L} -k_{R_L} e_{R_L} + [J_L e_{\Omega_L} +\label{Rerror} \\
    & (2J_L - \Tr(J_L) I)R_L^T \tilde{R}_L \tilde{\Omega}_L]^{\vee} e_{\Omega_L} - \sum_{j=1}^4 \hat{r}_j R_L^T (\tilde{q}_j^T \tilde{\mu}_j)\hat{q}_j e_{q_j},\nonumber\\
    &e_{\omega_j} = (\dot{q}_j^T \tilde{\omega}_j)q_j - k_{\omega} e_{\omega_j} - k_q e_{q_j}.\label{qerror}
\end{align}

\subsection{Theoretical Analysis}

The goal now is to prove that the gains can be chosen in such a way that the error dynamics \eqref{xerror}-\eqref{qerror} have an exponentially stable equilibrium point at the origin. This will be accomplished via Lyapunov analysis. Before this, however, we will introduce two configuration errors—$\Psi_{q_j}$ and $\Psi_{R_L}$. These are real-valued functions which are positive-definite about the points $q_j = \tilde{q}_j$ and $R_L = \tilde{R}_L$, respectively, and are defined explicitly as $\Psi_{q_j} = 1 - \tilde{q}_j^T q_j$ and $\Psi_{R_L} = \frac12\Tr\left(I - \tilde{R}_L^T R_L\right).$

We now define a Lyapunov candidate $V: \mathcal{D} \to \R$ defined as
\begin{align*}
    V =& \frac12 ||e_{v_L}||^2 + \frac12 k_{x_L} ||e_{x_L}||^2 + c_x e_{x_L}^T e_{v_L} \\
    & +\sum_{j=1}^4 \left[ \frac12 ||e_{\omega_j}||^2 + k_q \Psi_{q_j} + c_q e_{q_j}^T e_{\omega_j} \right]\\
    &+ \frac12 e_{\Omega_L}^T J_L \Omega_L + k_{R_L} \Psi_{R_L} + c_R e_{R_L}^T J_L e_{\Omega_L},
\end{align*} which is defined on the domain
\begin{align*}
\mathcal{D} &= \{(e_{x_L}, e_{v_L}, e_{R_L}, e_{\Omega_L}, e_{q_j}, e_{\omega_j}) \ : \ ||e_{x_L}|| < e_{x_{max}}, \\
&\, \qquad \Psi_{R_L} < \psi_{R_L} < 1, \ \Psi_{q_j} < \psi_{q_j} < 1 \}, 
\end{align*} where $c_x, c_q, c_R, \psi_{q_j}, \psi_{R_L}$ are positive real numbers.

Before establish sufficient conditions for the gains in order to achieve exponential stability of the zero equilibrium for the tracking error we  state the following Lemma that follows similarly as in \cite{jacobquads} and \cite{lee_geometric_2018}.

\begin{figure*}
\begin{align}
\underbar{P}_{x_L} &= \frac12 \begin{bmatrix} k_{x_L} & -c_x \\ -c_x & 1 \end{bmatrix}, \quad \bar{P}_{x_L} = \frac12 \begin{bmatrix} k_{x_L} & c_x \\ c_x & 1 \end{bmatrix},\,\quad \underbar{P}_{R_L} = \frac12 \begin{bmatrix} 2k_{R_L} & -c_{R} \lambda_{\max}(J_L) \\ -c_{R} \lambda_{\max}(J_L) & \lambda_{\min}(J_L) \end{bmatrix}, \label{matrices1}\\
\bar{P}_{R_L} &= \frac12 \begin{bmatrix} \frac{2k_{R_L}}{2 - \psi_{R_L}} & c_{R} \lambda_{\max}(J_L) \\ c_{R} \lambda_{\max}(J_L) & \lambda_{\max}(J_L) \end{bmatrix},\,
\underbar{P}_{q_j} = \frac12 \begin{bmatrix} 2k_q & -c_q \\ -c_q & 1 \end{bmatrix}, \bar{P}_{q_j} = \frac12 \begin{bmatrix} \frac{2k_q}{2 - \psi_{q_j}} & c_q \\ c_q & 1 \end{bmatrix}, W_{qR_j} = \begin{bmatrix} c_R B & 0 \\ \alpha_L \sigma_j k_{R_L} + B & 0 \end{bmatrix}\label{eqn:bM3}\\
W_j &= \begin{bmatrix} W_{x_j} & -\frac12 W_{xR_j} & -\frac12 W_{xq_j} \\ -\frac12 W_{xR_j} & W_{R_j} & -\frac12 W_{Rq_j} \\
-\frac12 W_{xq_j} & -\frac12 W_{Rq_j} & W_{q_j}\end{bmatrix},\,\quad
W_{x_j} = \frac14 \begin{bmatrix} c_x k_{x_L}(1 - 4\alpha_j \beta) & -\frac12 c_x k_{v_L} (1 + 4\alpha_j \beta) \\ -\frac12 c_x k_{v_L} (1 + 4\alpha_j \beta) & k_{v_L} (1 - 4\alpha_j \beta) - c_x \end{bmatrix},\\
W_{R_j} &= \frac14 \begin{bmatrix} c_R k_{R_L} (1 - 4\alpha_j \sigma_j) & -\frac{c_R}{2}(k_{\Omega_L} + B + 4\alpha_j \sigma_j) \\ -\frac{c_R}{2}(k_{\Omega_L} + B + 4\alpha_j \sigma_j) & k_{\Omega_L} (1 - 4 \alpha_j \sigma_j) - 2 c_R \lambda_{\max}(J_L)
\end{bmatrix},\,\quad W_{xq_j} = \begin{bmatrix} c_x B & 0 \\ \beta k_{x_L} e_{x_{\max}} + B & 0 \end{bmatrix}, \\
W_{q_j} &= \begin{bmatrix} c_q k_q & -\frac12 c_q (k_{\omega} + C_{q_j}) \\ -\frac12 c_q (k_{\omega} + C_{q_j}) & k_{\omega} - c_q \end{bmatrix},\quad W_{xR_j} = \alpha_j \begin{bmatrix} \gamma c_x k_{R_L} + \delta_j c_R k_{x_L} & \gamma c_x k_{\Omega_L} + \delta_j k_{x_L} \\ \gamma k_{R_L} + \delta_j c_R k_{v_L} & \gamma k_{\Omega_L} + \delta_j k_{v_L}
\end{bmatrix}.\label{matrices2}
\end{align}
\end{figure*}

\begin{lemma} Consider the set of matrices \eqref{matrices1}-\eqref{matrices2} for $\alpha_j := \sqrt{\psi_{q_j}(2 - \psi_{q_j})}$, \ $\alpha_L := \sqrt{\psi_{R_L}(2 - \psi_{R_L})}$, $\gamma = \frac1{m_L \gamma_{\min}(\mathcal{P}\mathcal{P}^T)}$, $\beta = m_L \gamma$, $\delta_j = m_L \frac{||\hat{r}_j||}{\sqrt{\lambda_{\min}(\mathcal{P}\mathcal{P}^T)}}, \sigma_j = \frac{\delta_j}{m_L}$, with $B$ a constant which can be obtained from $\tilde{x}_L$ and $\tilde{R}_L$.
If $\bar{P}_{x_L}, \underbar{P}_{x_L}, \bar{P}_{R_L}, \underbar{P}_{R_L}, \bar{P}_{q_j}, \underbar{P}_{q_j}$, and $\frac12(W_j + W_j^T)$ are positive-definite matrices for $j = 1,...,4$, then the origin of the error dynamics \eqref{xerror}-\eqref{qerror} is an exponentially stable equilibrium point.
\end{lemma}

Now we seek to find sufficient conditions under which a choice in gains exist such that these conditions are satisfied.
\begin{theorem}\label{Prop: gains}
Consider the control system defined by \eqref{xred}-\eqref{Rdynred} and control inputs determined by \eqref{cont1}-\eqref{cont3}, with $u_j^{\parallel}$ and $u_j^{\perp}$ given by equation \eqref{controllers}-\eqref{controllers2}. For sufficiently small $\alpha_j$, there exists gains $k_{x_L}$, $k_{v_L}$, $k_{R_L}$, $k_{\Omega_L}$, $k_q$, $k_{\omega}$, $j=1,\ldots,4$, such that the zero equilibrium of the tracking error $(e_{x_L}, e_{v_L}, e_{R_L}, e_{\Omega_L}, e_{q_j}, e_{\omega_j})$ is exponentially stable.
\end{theorem}
\textit{Proof:} It is easy to see that the matrices $\bar{P}_{x_L}, \underbar{P}_{x_L}, \bar{P}_{R_L}, \underbar{P}_{R_L}, \bar{P}_{q_j},$ and $\underbar{P}_{q_j}$ are positive-definite for $c_x, c_R, c_q$ sufficiently small. Now denote the symmetric part of $W_j$ by $\mathcal{W}_j = \frac12 (W_j + W_j^T)$ and similarly define the symmetric parts of the submatrices by $\mathcal{W}_{xR_j}, \mathcal{W}_{xq_j},$ and $\mathcal{W}_{qR_j}$. It is clear that $\mathcal{W}_j$ can be expressed in the form $\mathcal{W}_j= \begin{bmatrix} P & S \\ S^T & Q \end{bmatrix}$, where $P = \begin{bmatrix} W_{x_j} & -\frac12 \mathcal{W}_{xR_j} \\ -\frac12 \mathcal{W}_{xR_j} & W_{R_j} \end{bmatrix}$, $S = -\frac12 \begin{bmatrix} \mathcal{W}_{xq_j} & \mathcal{W}_{qR_j} \end{bmatrix}^T$, and $Q = {W}_{q_j}$. Now, observe that $\mathcal{W}$ can be decomposed as: \begin{align*}
\begin{bmatrix} P & S \\ S^T & Q \end{bmatrix} = \begin{bmatrix} I & SQ^{-1} \\ 0 & I \end{bmatrix} \begin{bmatrix} P - SQ^{-1}S^T  & 0 \\ 0 & Q \end{bmatrix} \begin{bmatrix} I & SQ^{-1} \\ 0 & I \end{bmatrix}^T
\end{align*}where $P - SQ^{-1}S^T$ is often referred to as the \textit{Schur complement} of $Q$. It then follows that $\mathcal{\mathcal{W}} \succ 0$ if and only if $P - SQ^{-1} S^T \succ 0$ and $Q \succ 0.$ Note that $P - SQ^{-1}S^T$ can itself be expressed in form of a $4\times4$ block matrix given by \begin{align*}
    &P - SQ^{-1}S^T = \\
&{\footnotesize\begin{bmatrix}
W_{x_j} - \frac14 \mathcal{W}_{xq_j} {W}_{q_j}^{-1} \mathcal{W}_{xq_j} & -\frac12 \mathcal{W}_{xR_j} - \frac14 \mathcal{W}_{qR_j}{W}_{q_j}^{-1}\mathcal{W}_{xq_j} \\ -\frac12 \mathcal{W}_{xR_j} - \frac14 \mathcal{W}_{qR_j}{W}_{q_j}^{-1}\mathcal{W}_{xq_j} & W_{R_j} -\frac14 \mathcal{W}_{qR_j}{W}_{q_j}^{-1}\mathcal{W}_{qR_j}
\end{bmatrix}}.
\end{align*}
Repeating the previous analysis, but now on $P - SQ^{-1}S^T$, we find that $\mathcal{W}_j\succ 0$ if and only if the following three conditions hold: 
$(\bm{1}) W_{q_j} \succ 0$, $(\bm{2}) W_{R_j} - \frac14 \mathcal{W}_{qR_j} W_{q_j}^{-1} \mathcal{W}_{qR_j} \succ 0$ and $(\bm{3})$\footnotesize{\begin{align*}
0&\prec W_{x_j} - \frac14 \mathcal{W}_{xq_j} W_{q_j}^{-1} \mathcal{W}_{xq_j} \\&- \frac14 \mathcal{W}_{xR_j} (W_{R_j} - \frac14 \mathcal{W}_{qR_j} W_{q_j}^{-1} \mathcal{W}_{qR_j})^{-1} \mathcal{W}_{xR_j} \\
& - \frac18 \mathcal{W}_{xR_j} (W_{R_j} - \frac14 \mathcal{W}_{qR_j} W_{q_j}^{-1} \mathcal{W}_{qR_j})^{-1} \mathcal{W}_{qR_j}W_{q_j}^{-1}\mathcal{W}_{xq_j} \\
& - \frac18 \mathcal{W}_{xq_j} W_{q_j}^{-1} \mathcal{W}_{qR_j} (W_{R_j} - \frac14 \mathcal{W}_{qR_j} W_{q_j}^{-1} \mathcal{W}_{qR_j})^{-1} \mathcal{W}_{xR_j} \\
& - \frac1{16} \mathcal{W}_{xq_j} W_{q_j}^{-1} \mathcal{W}_{qR_j} (W_{R_j} - \frac14 \mathcal{W}_{qR_j} W_{q_j}^{-1} \mathcal{W}_{qR_j})^{-1} \mathcal{W}_{qR_j}W_{q_j}^{-1}\mathcal{W}_{xq_j}.
\end{align*}}
 
\normalsize Moreover, the minimum (maximum) eigenvalue of $\mathcal{W}$ is exactly the smallest (largest) of the minimum (maximum) eigenvalues of the three matrices in the above conditions. We now seek to verify that appropriate choices in the gains and constants can be made to satisfy the above conditions. First, by looking at the characteristic equation of $W_{q_j}$, 
\begin{align*}
2\lambda_{\min}(W_{q_j}) =& (k_{\omega} - c_q + c_qk_q) \\
&- \sqrt{(k_{\omega} - c_q - c_qk_q) + c_q^2(k_{\omega} + C_{q_j})^2 }.
\end{align*}

By taking $c_q$ sufficiently small, $k_{\omega}$ sufficiently large, and $k_q = \frac{k_{\omega}}{c_q}$, $\lambda_{\min}(W_{q_j})$ can be made arbitrarily large. Consequently, $\lambda_{\max}(W_{q_j}^{-1}) = \lambda_{\min}(W_{q_j})^{-1}$ can be made arbitrarily small (and positive). Hence we have that $\frac14 \mathcal{W}_{qR_j} W_{q_j}^{-1} \mathcal{W}_{qR_j} \succeq 0$ and, since $\mathcal{W}_{qR_j}$ and $W_{q_j}$ are independent, we can shrink the maximum eigenvalue of $\mathcal{W}_{qR_j} W_{q_j}^{-1} \mathcal{W}_{qR_j}$ arbitrarily by shrinking the maximum eigenvalue of $W_{q_j}^{-1}$.

Similarly, from the characteristic equation of $W_{R_j}$, we see that the eigenvalues satisfy
\begin{align*}&2\lambda_{\min} = (c_R k_{R_L} + k_{\Omega_L})\nu_j - 2c_R \lambda_{\max}(J_L) - \\
& \sqrt{\left((c_R k_{R_L} - k_{\Omega_L})\nu_j + 2 c_R \lambda_{\max}(J_L) \right)^2 + c_R^2\left(k_{\Omega_L} + B + \nu_j\right)^2}
\end{align*} 
where $\nu_j = 1 - 4\alpha_j \sigma_j > 0$ for $\alpha_j$ sufficiently small. Now choose $k_{\Omega_L} = c_R k_{R_L} + \frac{2 c_R \lambda_{\max}(J_L)}{\nu_j}$. Then,
\begin{align*}
    \lambda_{\min} &= c_R \left[ 2k_{R_L} (\nu_j - \frac12 c_R) - B - 4\alpha_j \sigma_j - \frac{2c_R \lambda_{\max}(J_L)}{\nu_j}\right].
\end{align*}
Therefore, if we choose $c_R$ sufficiently small so that $\nu_j - \frac12 c_R > 0$, it follows that the minimum eigenvalue of $W_{R_L}$ can be arbitrarily large by taking $k_{R_L}$ sufficiently large. It then follows that $W_{R_j} - \frac14 \mathcal{W}_{qR_j} W_{q_j}^{-1} \mathcal{W}_{qR_j} \succ 0$, and its minimum eigenvalue can be made arbitrarily large with appropriate choices of $k_{R_j}, k_q, k_{\omega}$ for $j = 1,...,4$.

Now, we look at condition $(\bm{3})$. First, choose $k_{x_L}, k_{v_L}, c_x$ such that $W_{x_j} \succ 0$ (this can always be done by, for instance, choosing $c_x$ sufficiently small). We now wish to show that the remaining subtractive terms can be shrunk arbitrarily independent of $W_{x_j}$. Observe that $\mathcal{W}_{xq_j}$ and $W_{q_j}^{-1}$ are independent, so that we may force the maximum eigenvalue of $\mathcal{W}_{xq_j} W_{q_j}^{-1} \mathcal{W}_{xq_j} \succeq 0$ to be arbitrarily small after (potentially) further shrinking the maximum eigenvalue of $W_{q_j}^{-1}$.

Further observe that we may write the third term, $\frac14 \mathcal{W}_{xR_j} (W_{R_j} - \frac14 \mathcal{W}_{qR_j} W_{q_j}^{-1} \mathcal{W}_{qR_j})^{-1} \mathcal{W}_{xR_j}$, in the form $\alpha^2 M^T A M$, where $M$ is independent of $\alpha$ and the terms of $A$ are at most of $\mathcal{O}(\frac1{\alpha})$. Hence, we may shrink this term arbitrarily by shrinking $\alpha$. 

The fourth and fifth terms are transposes of each other and therefore may be handled simultaneously. Note that find that the maximum eigenvalue is bounded above by
\footnotesize{\begin{align*}
    \lambda_{\max}&((W_{R_j} - \frac14 \mathcal{W}_{qR_j} \mathcal{W}_{q_j}^{-1}\mathcal{W}_{qR_j})^{-1})||\mathcal{W}_{qR_j} \mathcal{W}_{q_j}^{-1} \mathcal{W}_{xq_j}|| \ ||\mathcal{W}_{xR_j}|| \\
    &\le \sqrt{\lambda_{\max}(\mathcal{W}_{q_j}^{-1})} \lambda_{\max}((W_{R_j} - \frac14 \mathcal{W}_{qR_j} \mathcal{W}_{q_j}^{-1}\mathcal{W}_{qR_j})^{-1})\\&\times ||\mathcal{W}_{qR_j}|| \ || \mathcal{W}_{xq_j}|| \ ||\mathcal{W}_{xR_j}||.
\end{align*}}
\normalsize This term therefore can be arbitrarily shrunk by shrinking the maximum eigenvalue of $\mathcal{W}_{q_j}$. Moreover, the presence of the norm $||W_{xR_j}||$ also gives us control of the size of the term via $\alpha.$ The final term is handled similarly—we find that the maximum eigenvalue is bounded  by \footnotesize{\begin{align*}
   \lambda_{\max}&((W_{R_j} - \frac14 \mathcal{W}_{qR_j} \mathcal{W}_{q_j}^{-1}\mathcal{W}_{qR_j})^{-1})||\mathcal{W}_{qR_j} \mathcal{W}_{q_j}^{-1} \mathcal{W}_{xq_j}||^2 \\
    \le &\lambda_{\max}(\mathcal{W}_{q_j}^{-1}) \lambda_{\max}((W_{R_j} - \frac14 \mathcal{W}_{qR_j} \mathcal{W}_{q_j}^{-1}\mathcal{W}_{qR_j})^{-1})\\&\times ||\mathcal{W}_{qR_j}||^2 \ || \mathcal{W}_{xq_j}||^2,
\end{align*}}\normalsize which again may be shrunk arbitrarily by shrinking the maximum eigenvalue of $\mathcal{W}_{q_j}^{-1}$. So, for sufficiently small $\alpha$ and $\lambda_{\max}(W_{q_j})$ sufficiently large, condition $(\bm{3})$ holds and $\mathcal{W}_j \succ 0$. 

Now denote by $S^5$ the unit sphere in $\mathbb{R}^6$ (i.e., $x\in\mathbb{R}^6$ such that $||x||=1)$, and suppose that $x \in S^5$ is such that $x^T \mathcal{W}_jx = \lambda_{\min}(\mathcal{W}_j)$, and decompose $x = (x_1, x_2) \in \R[4] \times \R[2]$. Then, from the Schur decomposition, we have $\lambda_{\min}(\mathcal{W}) \ge \min\{ \lambda_{\min}(Q), \ \lambda_{\min}(P - SQ^{-1}S^T)\} \big|\big|x + \begin{bmatrix} Q^{-1}S^Tx_1 & 0 \end{bmatrix}^T\big|\big|^2$. Note that $||Q^{-1}S^Tx_1||$ can be made arbitrarily small by increasing the maximum eigenvalue of $Q$, so that $\big|\big|x + \begin{bmatrix} Q^{-1}S^Tx_1 & 0 \end{bmatrix}^T\big|\big|$ can be made arbitrarily close to $1$. Repeating this procedure with $P - SQ^{-1}S^T$, we find that the minimum eigenvalue of $\mathcal{W}_j$ is bounded below by a quantity that can be made arbitrarily close to the minimum of the minimum eigenvalue of the matrices in conditions $(\bm{1})$, $(\bm{2})$, and $(\bm{3})$ above—all of which can be made arbitrarily large. Hence, the minimum eigenvalue of $\mathcal{W}_j$ can be made arbitrarily large.
$\hfill\square$

\subsection{Control design for the full model}

We now wish to include the quadrotor attitude dynamics back into the problem, so that we have a control scheme for the full reduced model. Note however that in the design of the geometric controllers $u_j$, we assumed that each quadrotor can generate a thrust along any direction. However, the quadrotor dynamics are in fact underactuated, since the direction of the total thrust is always parallel to its third body-fixed axis. Hence we desire the attitude of each quadrotor to be controlled such that the third body-fixed axis becomes parallel to the direction of the control force $u_j$. This was accomplished for instance in \cite{lee_geometric_2010}, yielding the moment controller
\begin{align}
M_j = &-\frac{k_{R_j}}{\epsilon^2} e_{R_j} -\frac{k_{\omega}}{\epsilon} e_{\Omega_j} + \Omega_j\times J_j\Omega_j\nonumber\\
& -J_j (\hat\Omega_j R_j^T \tilde{R}_j\tilde{\Omega}_j-R_j^T \tilde{R}_j\dot{\tilde{\Omega}}_j),\label{eqn:Mi}
\end{align} for $j=1,...,4$ and where $\epsilon,k_{R_j},k_{\omega}$ are positive constants. 



Stability of the corresponding controlled system for the unreduced model can be studied by using singular perturbation theory for the attitude dynamics of quadrotors as in~\cite{sreenath_geometric_2013}, \cite{lee_geometric_2013}. In particular, as a direct application of \cite{lee_geometric_2010}, \cite{lee_geometric_2018}, in the context of Theorem \ref{Prop: gains} for our particular cooperative transportation task lead to the following result.

\begin{corollary}\label{corol}
Consider the control system defined
by \eqref{xred}-\eqref{Rdynred}  with control inputs determined by \eqref{cont1}-\eqref{cont3}, with $u_j^{\parallel}$ and $u_j^{\perp}$ given by equation \eqref{controllers}-\eqref{controllers2}. There
exists $\delta > 0$, such that for all $\epsilon< \delta$, the zero equilibrium of the tracking errors $(e_{x_L}, e_{v_L}, e_{R_L}, e_{\Omega_L}, e_{q_j}, e_{\omega_j}, e_{R_j}, e_{\Omega_j})$ is exponentially stable.
\end{corollary}

Now that we have established the exponential tracking of the full reduced model, we wish to connect this back to our original model with elastic cables. This can be done by showing that our system is under the conditions of Theorem 11.2 in \cite{khalil}. Before stating the Proposition formally, we introduce some notation and definitions needed to introduce the formal result. 

\begin{definition}
The \textit{boundary layer system} for the singular perturbation problem given by \eqref{spmodel} is defined as:
$$\frac{\partial r}{\partial \tau} = g(t,x,r+h(t,x), 0),$$where $r := z - h(t,x)$ with $h(t,x)$ as defined by \eqref{eqh} and $\tau := \frac{t - t_0}{\epsilon}$ for $t_0$ the value of time from which we obtain our initial data. 
\end{definition}

The following Corollary for the exponential stability of the boundary layer system \eqref{spmodel} follows from the case of a single quadrotor transporting a point mass load with an elastic cable (see Lemma 2 in \cite{sreenath_geometric_2017}). 

\begin{corollary}\label{boundarylayer}
The boundary layer system for \eqref{spmodel} with control inputs $u_j$ and $M_j$ as defined above has an exponentially stable equilibrium point at the origin.
\end{corollary}

Theorem 11.2 in \cite{khalil} tells us that the trajectories of the original model lie in a neighborhood of the trajectories of the reduced model when the origin of the boundary layer system and the error dynamics of the reduced model are exponentially stable -- which follows immediately from Corollary \ref{boundarylayer} and Theorem \ref{Prop: gains} above. Formally stated, we have the following:

\begin{proposition}\label{singular}
Let the control inputs $u_j$ and $M_j$ be defined as above. Denote by $x(t)$ a trajectory of the reduced model \eqref{xred}-\eqref{Rdynred} which converges exponentially to the desired trajectory. Denote by $r(t)$ a trajectory of the boundary layer system which converges exponentially to the origin. 

There exists a positive constant $\epsilon^\ast$ such that for all $t \ge t_0$ and $0 < \epsilon < \epsilon^\ast$, there exists a unique solution $x(t, \epsilon)$, $z(t, \epsilon)$ of the singular perturbation problem \eqref{spmodel} on $[t_0, \infty)$ satisfying $
x(t, \epsilon) - x(t) = \mathcal{O}(\epsilon)$ and 
$z(t, \epsilon) - h(t, x(t)) - r\left(\frac{t - t_0}{\epsilon}\right) = \mathcal{O}(\epsilon)$, uniformly on $t \in [t_0, \infty)$. Moreover, for $t_1 > t_0$, we have  $z(t,\epsilon) - h(t,x(t)) = \mathcal{O}(\epsilon)$, uniformly on $(t_1, \infty)$ for $\epsilon < \epsilon^{\ast\ast} < \epsilon^\ast$.
\end{proposition}
\subsection{Control design in the presence of unstructured disturbances}\label{sec_dist}
It is well-known that exponential stability is robust to small disturbances \cite{khalil2}, and therefore the control scheme developed in Section \ref{sec5} will hold provided that the external disturbances and error measurements are sufficiently small. This can not always be guaranteed in real life applications, so to that end we will introduce bounded but unstructured disturbances to the reduced model \eqref{xred}-\eqref{wjred} and test the previous control scheme in this new scenario. 

In particular, for $j=1,\ldots,4$, we consider the following perturbed dynamical system
\begin{align}
&\dot{x}_L = v_L,\,\, \dot{R}_L = R_L \hat{\Omega}_L,\,\, \dot{q}_j = \omega_j \times q_j,\label{kinpert}\\ 
&m_{\eff} ( \dot{v}_L + ge_3) = \\ &\hspace{1cm}\sum_{j=1}^4 \Big(u_j - m_Q R_L(\hat{\Omega}_L^2 + \dot{\hat{\Omega}}_L)r_j + m_Q L\ddot{q}_j \Big) + \Delta x_L, \nonumber\\
&J_{\eff} \dot{\Omega}_L + \hat{\Omega}_L J_{\eff} \Omega_L =\nonumber\\ &\hspace{.1cm} \sum_{j=1}^4 m_Q \hat{r}_j R_L^T \Big(-g e_3 - \dot{v}_L + L\ddot{q}_j + \frac1{m_Q}u_j\Big)+ \Delta R_L,\label{Omegadist} \\ 
&\dot{\omega}_j = L^{-1} \hat{q}_j\Big( \dot{v}_L - R_L(\hat{\Omega}_L^2 + \dot{\hat{\Omega}}_L)r_j + ge_3 - m_Q^{-1} u_j \Big)+ \Delta q_j\label{wjdist}
\end{align} 
where $\Delta x_L, \Delta R_L, \Delta q_j$ are unstructured disturbances satisfying the bounds $||\Delta x_L|| \le \bar{x}_L$, $||\Delta R_L|| \le \bar{R}_L$ and $||\Delta q_j|| \le \bar{q}_j$ for some real numbers $\bar{x}_L, \bar{R}_L, \bar{q}_j$.

Defining the controllers, configuration errors, and Lyapunov candidate $V$ as in Section \ref{sec5}, it is easy to see that the matrices $\bar{P}_{x_L}, \underbar{P}_{x_L}, \bar{P}_{R_L}, \underbar{P}_{R_L}, \bar{P}_{q_j}, \underbar{P}_{q_j}$ remain unchanged. On the other hand, the upper bound on the time derivative of $V$ is modified as $\dot{V} \le -\sum_{j=1}^4 z_j^T W_jz_j + E^T z,$ where $E := \begin{bmatrix}
\frac{c_x\delta x_L}{m_r}, \frac{\delta x_L}{m_r}, \frac{3c_R \delta R_L}{2m_r L_r}, \frac{3\delta R_L}{2m_r L_r}, \frac{c_{q_i}\delta q_j}{m_Q L_c}, \frac{\delta q_j}{m_Q L_c} 
\end{bmatrix}$, \ $z_j = \begin{bmatrix} ||e_{x_L}|| & ||e_{v_L}|| & ||e_{R_L}|| & ||e_{\Omega_L}|| & ||e_{q_j}|| & ||e_{\omega_j}||   \end{bmatrix}$ and $W_j$ is as in Section \ref{sec5}. Fix $\epsilon > 0.$ 
From Young's Inequality, we have $E^T z \le \frac{||E||^2}{16 \epsilon} + 4\epsilon ||z||^2$, and hence
$$\dot{V} \le -\sum_{j=1}^4 z_j^T (W_j- {\epsilon} I) z_j + \frac{||E||^2}{16 \epsilon}.$$
Analogous to before, we replace the matrix $W_j- \epsilon I$ with its symmetric part $\mathcal{W}^\ast_j := \frac12(W_j+ W_j^T) - {\epsilon} I$. Note that $\lambda_{\min}(\mathcal{W}^\ast) = \min_{j=1,...,4}\lambda_{\min}(\mathcal{W}_j) - \epsilon$. From Corollary \ref{Prop: gains}, $\lambda_{\min}(\mathcal{W})$ can be made arbitrarily large by choosing gains appropriately. Hence, we can choose them so that $\mathcal{W}^\ast_j \succ 0$, with an arbitrarily large minimum eigenvalue. We then have 
\begin{align*}
\lambda_{\min}(\underbar{P})||z||^2 &\le V \le \lambda_{\max}(\bar{P})||z||^2, \\
\dot{V} &\le -4\lambda_{\min}(\mathcal{W}^\ast_j)||z||^2 + \frac{||E||^2}{16\epsilon}.
\end{align*}
This implies that $\dot{V} \le -\frac{4\lambda_{\min}(\mathcal{W}^\ast_j)}{\lambda_{\max}(\bar{P})}V + \frac{||E||^2}{16\epsilon}$, so that $\dot{V} < 0$ when $V > \frac{\lambda_{\max}(\bar{P})}{\lambda_{\min}(\mathcal{W}^\ast_j)}\frac{||E||^2}{64\epsilon}:= d_1 > 0.$ Clearly, $d_1$ can be shrunk arbitrarily by increasing $\lambda_{\min}(\mathcal{W}^\ast_j)$. If we now define the set $S_{r} := \{z \in \mathcal{D} : V(z) < r \}$, where $r$ is some real number, then any trajectory starting in the open set $\mathcal{D} \setminus \bar{S}_{d_1}$ will converge exponentially to the region $\bar{S}_{d_1}$, where $\bar{S}_{d_1}$ denotes the topological closure of $S_{d_1}$. Since $V$ is continuous and positive, $\bar{S}_{d_1}$ is some closed neighborhood of the origin that can be made arbitrarily small (by shrinking $d_1$). We formalize this result with the following theorem:

\begin{theorem}\label{gainsthdist}
Consider the system with disturbances defined by equations \eqref{kinpert}-\eqref{wjdist} with control inputs determined by \eqref{cont1}-\eqref{cont3}, with $u_j^{\parallel}$ and $u_j^{\perp}$ given by equation \eqref{controllers}-\eqref{controllers2}. For sufficiently small $\alpha$, there exists control gains $k_{x_L}$, $k_{v_L}$, $k_{R_L}$, $k_{\Omega_L}$, $k_q$, $k_{\omega}$, $j=1,\ldots,4$, such that the zero equilibrium of the tracking errors $e_{x_L}$, $e_{v_L}$, $e_{R_L}$, $e_{\Omega_L}$, $e_{q_j}$, $e_{\omega_j}$, are uniformly ultimately bounded with an arbitrarily small ultimate bound. 
\end{theorem}

\begin{remark}
Note that, as we did in the previous subsection, the next step must be include the attitude and dynamics for quadrotors and back to the elastic model, nevertheless, we will not extend this result to include the quadrotor attitude kinematics and dynamics here, nor relate it back to the original elastic model since such a task follows essentially the same strategy as it did in Section \ref{sec5}, just replacing exponential stability with uniform ultimate boundedness where appropriate.  
\end{remark}
\section{Numerical Validation}

To validate the controller we first consider a numerical simulation with four quadrotors carrying a load with a mass of $\SI{2.0}{\kilogram}$, and rigid body shape of rectangular box with length, width and height are $\SI{2.0}{\meter}$, $\SI{1.0}{\meter}$ of $\SI{0.2}{\meter}$, respectively. The inertial tensor of the load can be computed due to its shape and it is defined by $J_L = \diag\left[1.04, 5.0, 4.04\right]^T$. The mass of each quadrotor is $m_Q =\SI{0.755}{\kilogram}$ The length of the cables are set to $l_j=\SI{1.0}{\meter}$ and the points where the load was attached are given by $r_1=\left[0.5, 1.0, 0.1\right], r_2=\left[0.5, -1.0, 0.1\right], r_3=\left[-0.5, -1.0, 0.1\right]$ and $r_4=\left[-0.5, 1.0, 0.1\right]$, where each point described before is referred from the center of the mass of the load to the center of mass of each quadrotor. The numerical simulation was realized with $\SI{40}{\second}$ of simulation, with time step of $\SI{2}{\milli\second}$ implemented in Python on a computer with a processor Intel Core i7-6800K at $\SI{3.4}{\giga\hertz}$ and $\SI{64}{\giga\byte}$ of RAM memory. The desired trajectory for position of the load used for this simulation is given by $
x_d (t) = \left[1.2 \sin(0.4 \pi t), 4.2 \cos(0.2 \pi t), 5.0 \right]$, and the desired attitude by $\displaystyle{R_d (t) = \left[\frac{\tilde{v}}{|\tilde{v}|}, \frac{\hat{e_3}\tilde{v}}{|\hat{e_3}\tilde{v}|}, e_3 \right]}$, where $\tilde{v}$ is the desired velocity. Initial conditions are given by 
$x(0) = \left[1.5, 2.5, 2.5 \right]^T\si{\meter}$, $R(0) = \diag(1,1,1)$ together with linear and angular velocities starting at rest.

\begin{figure}[htb]
\centering
\includegraphics[width=8cm]{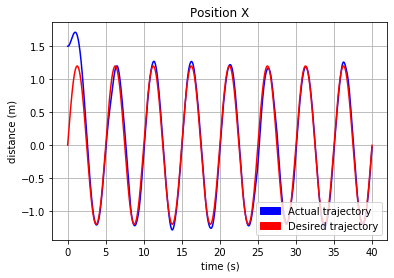}\\
\includegraphics[width=8cm]{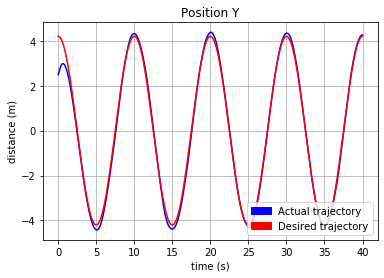}\\
\includegraphics[width=8cm]{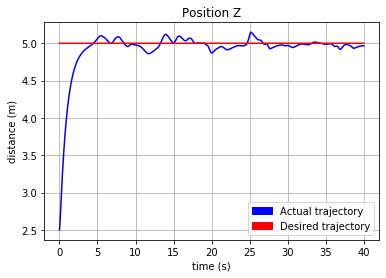}
\caption{Response of controlled system. In top shows X response, on middle Y response and on bottom Z response for the position of the load.} 
\label{fig:non_pos}
\end{figure}

The control gains are set to $k_{x_L}=k_{v_L}=600$. This values for control gains were selected to guarantee the global exponentially stability for the system by using conditions $(\bm{1})$, $(\bm{2})$ and $(\bm{3})$ in the proof of Theorem \ref{Prop: gains}.  The response of the controlled system is visualized in Figure \ref{fig:non_pos}. 

\begin{figure}[htb]
\centering
\includegraphics[width=8.3cm]{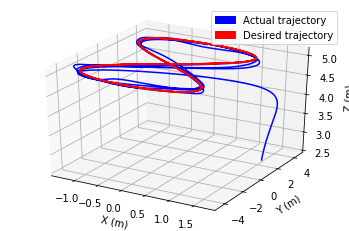}
\includegraphics[width=8.3cm]{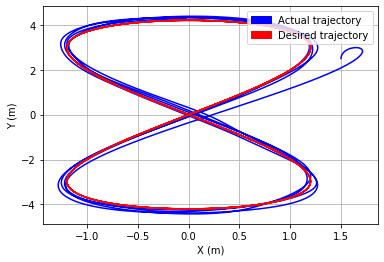}
\caption{Trajectory tracking for load position.} 
\label{fig:non_trajectory}
\end{figure}

The Mean Squared Error (MSE) for each axis are $e_{X_{MSE}}=0.036, e_{Y_{MSE}}=0.051$ and $e_{Z_{MSE}}=0.086$. The trajectory tracking executed by the controller is shown in Figure \ref{fig:non_trajectory}.

Next, we analyze the control performance in a simulation where the system is under disturbances. The disturbance, as introduced in Section \ref{sec_dist}, is set to

\begin{align*}
    \Delta x_L (t) &= \begin{bmatrix}
    0.5 \sin(0.43t)\\
    0.5 \cos(0.21t)\\
    0.2 \sin(0.75t)-e^{-t}\\
    \end{bmatrix},\\  \Delta R_L (t) &= \begin{bmatrix}
    0.2+0.45 \sin(3.0t)\\
    0.3 - 0.65 \cos(1.4t)\\
    0.05\sin(2.1t)\\
    \end{bmatrix}
\end{align*}

\begin{figure}[htb]
\centering
\includegraphics[width=8.5cm]{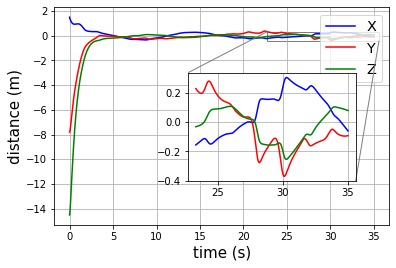}
\caption{Errors obtained from the numerical validation under disturbances $\Delta x_L$ and $\Delta R_L$.} 
\label{fig:error_pos}
\end{figure}
\begin{figure}[htb]
\centering
\includegraphics[width=8.5cm]{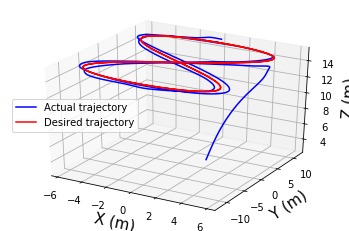}
\includegraphics[width=8.5cm]{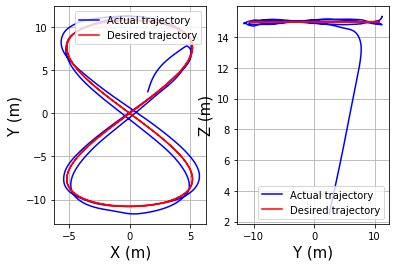}
\caption{Trajectory generated by the system under disturbances. On top, 3D visualization and bottom two views top and lateral.} 
\label{fig:trajectory}
\end{figure}

The tracking error is visualized in Figure \ref{fig:error_pos}, where we can observe the boundedness of the tracking error as stated in Theorem \ref{gainsthdist}. 
The MSE of the tracking error is given by $e_{X_{MSE}}=0.0398, e_{Y_{MSE}}=0.0804$ and $e_{Z_{MSE}}=0.1048$. Doing a comparative with the values obtained in the case without disturbances, the error increased due to the disturbances, but the controller still guarantees the boundedness of the tracking error. The trajectory made by the system under disturbance is shown in Figure \ref{fig:trajectory}.

\section{Conclusions}
We propose a geometric trajectory tracking controller for the cooperative task of four quadrotor UAVs transporting a rigid body load via inflexible elastic cables. This is handled in three stages: (i) Reduction of the model to that of a similar model with inelastic cables. We accomplish this by assuming sufficient stiffness and damping of the cables and utilizing the results of singular perturbation theory. (ii) Analysis of a geometric tracking controller in the reduced model. Lyapunov analysis is used to find sufficient conditions for stability, and Theorem \ref{Prop: gains} proves the existence of gains satisfying these conditions for sufficiently small initial errors in the cable attitudes. (iii) Under the same control law—trajectories of the original (elastic) model converge uniformly to the trajectories of the reduced model as the stiffness and damping of the cables approach infinity. (iv) Finally, we also extended the proposed approach to design a control law that guarantees bounded of the tracking error under unstructured bounded disturbances. 

In our model cables are attached to the center of each quadrotor. It would be interesting to explore in future work how to shift those attachment points and study how to deal with the resulting coupled systems—instead of a decoupled system. In addition, we are currently working to add uncertainties in order to further explore the robustness of the proposed controller. We also plan to study the construction of force variational integrators in optimal control problems, in a similar fashion to \cite{colombo2016geometric} and \cite{colombo2015variational}, dynamic interpolation problems \cite{bloch2021dynamic}, and obstacle avoidance problems \cite{bloch2017variational} for the cooperative task between quadrotors UAVs presented in this paper. 

\section*{Acknowledgments}
The project that gave rise to these results received the support of a fellowship from ”la Caixa” Foundation (ID 100010434). The fellowship codes are LCF/BQ/PI19/11690016 and LCF/BQ/DI19/11730028. The authors also acknowledge financial support from the Spanish Ministry of Science and Innovation, under grants PID2019-106715GB-C21. All the results of the paper are original and have not been presented nor submitted to a conference. The authors also acknowledge Manuela Gamonal from ICMAT, Spain, for the help provided with the numerical simulations.
\appendices

\section*{Appendix}

\subsection{Proof of Proposition 2.1
}
\textit{Proof:} We wish to apply Lagrange-d'Alembert Variational Principle for the Lagrangian $\mathbf{L}$ and external forces. Therefore, our system dynamics must satisfy

\begin{align}
&\delta \int_{0}^{T} \mathbf{L}(c(t),\dot{c}(t))\, dt \label{action}\\ &\hspace{0.9cm}+ \sum_{j=1}^{4} \int_{0}^{T} \left(\delta x_{Q_j}^T u_j + \left<R_j^T \delta R_j, \hat{M}_j\right> - c \dot{l}_j\delta l_j\right)dt = 0\nonumber
\end{align} where the integral on the right represents the virtual work done by the thrust controls $u_j$, the moment controls $M_j \in \R[3]$, and the spring damping, respectively. 

Expanding the variations within \eqref{action}, substituting the corresponding infinitesimal variations, and grouping like terms, we obtain
\begin{align*}
0&= \int_0^{T} \delta \dot{x}_L^T (m_{\eff} \dot{x}_L - m_Q \sum_{j=1}^4(\dot{\zeta}_j  - \dot{R}_L r_j))\\&+\delta x_L^T (-m_{\eff} ge_3 +\sum_{j=1}^{4}u_j)\,dt \\
&+ \  \int_0^{T} \sum_{j=1}^4 m_Q \delta \dot{R}_L r_j^T (\dot{x}_L + \dot{R}_L r_j - \dot{\zeta}_j) + \delta \Omega_L^T J_L.\\ &- m_Q g e_3^T (\delta R_L r_j) + \sum_{j=1}^4 u_j^T\delta R_L r_j  \\
&- \sum_{j=1}^{4}\int_0^{T} \left[ m_Q(\delta \dot{l}_j) q_j^T \dot{x}_{Q_j} + \delta l_j ( m_Q \dot{q}_j^T \dot{x}_{Q_j}\right.\\& \left.- m_Q ge_3^T q_j  - k(L - l_j) + c \dot{l}_j + q_j^T u_j ) \right] dt \\
&- \sum_{j=1}^{4}\int_0^{T} \left[\xi_j^T \left( q_j \times (m_Q \dot{l}_j \dot{x}_{Q_j} - m_Q g l_j e_3 + l_j u_j)\right.\right.\\&\left.\left. + \dot{q}_j \times m_Q l_j \dot{x}_{Q_j} \right) + \dot{\xi}_j^T (q_j \times m_Q l_j \dot{x}_{Q_j}) \right] dt \\
&+\sum_{j=1}^{4} \int_0^{T}  \left[ \dot{\eta}_j^T J_Q \Omega_j  + \eta_j^T \left(  J_Q \Omega_j \times \Omega_j + M_j \right) \right]dt,
\end{align*}where $m_{\eff} := 4m_Q + m_L, \ J_{\eff} := J_L - \sum_{j=1}^4 m_Q \hat{r}_j^2$, and $\zeta_{j} := l_j q_j$. 

Integrating by parts and applying the equality of mixed partial derivatives and the fact that variations vanish on the endpoints, we obtain
 \begin{align*}
0&= \int_{0}^{T} \delta x_L \left[ m_{\eff} (\ddot{x}_L + ge_3) -\sum_{j=1}^4 m_Q (\ddot{\zeta}_j + \dot{R}_L r_j)-u_j \right] dt \\&+ \sum_{j=1}^{4}\int_0^{T} \eta_j^T \left[ J_Q \dot{\Omega}_j - J_Q \Omega_j \times \Omega_j - M_j \right]dt \\
&+ \  \int_0^T \eta_L^T \left[  J_{\eff}(\Omega_L \times \ddot{\Omega}_L) - \Omega_L \times \left(m_Q\hat{r}_j(\sum_{j=1}^{4}\ddot{\zeta}_{j} + u_j\right.\right.\\&\left.\left.\hspace{2cm}-ge_3-\ddot{x}_L)\right)\right]\,dt \\
&- \sum_{j=1}^{4}\int_0^{T} \delta l_j \left[ m_Q q_j^T( \ddot{x}_L + R_L(\hat{\Omega}_L^2 + \dot{\hat{\Omega}}_L)r_j - \ddot{\zeta}_j+ ge_3) \right.\\& \left. - c \dot{l}_j+ k(L - l_j) - q_j^T u_j \right] dt \\
&- l_j \sum_{j=1}^{4}\int_0^{T} \xi_j^T \left[m_Q q_j \times ( \ddot{x}_L + R_L(\hat{\Omega}_L^2 + \dot{\hat{\Omega}}_L)r_j - \ddot{\zeta}_j + ge_3)\right.\\&\left. - q_j \times u_j \right]dt.
\end{align*}

Each of these integrals can be treated independently, as their respective variations are independent. That is, for the above equation to be satisfied, we necessarily have that each integral vanish identically. Applying the Fundamental Lemma of the Calculus of Variations \cite{GF} to each integral yields the dynamical system:
\begin{align*}
&m_{\eff} ( \dot{v}_L + ge_3) = \sum_{j=1}^4 u_j + m_Q \ddot{\zeta}_j - m_Q R_L(\hat{\Omega}_L^2 + \dot{\hat{\Omega}}_L)r_j\\ 
 &J_{\eff} \dot{\Omega}_L + \hat{\Omega}_L J_{\eff} \Omega_L = \sum_{j=1}^4 m_Q \hat{r}_j R_L^T (-g e_3 - \dot{v}_L \\& \hspace{5.9cm}+ \ddot{\zeta}_j + \frac1{m_Q}u_j), \nonumber \\ 
&m_Q q_j^T \ddot{\zeta}_j = m_Q q_j^T ( \dot{v}_L + R_L(\hat{\Omega}_L^2 + \dot{\hat{\Omega}}_L)r_j + ge_3 - \frac1{m_Q}u_j) \nonumber\\& - c \dot{l}_j  + k(L - l_j), \label{ldyn}\\ 
&q_j \times \ddot{\zeta}_j = q_j \times ( \dot{v}_L + R_L(\hat{\Omega}_L^2 + \dot{\hat{\Omega}}_L)r_j + ge_3 -\frac1{m_Q}u_j) \\ 
&J_Q \dot{\Omega}_j = J_Q \Omega_j \times \Omega_j + M_j,\,\,\dot{R}_j = R_j \hat{\Omega}_j,\,\, \ j = 1,\ldots,4,
\end{align*} where we have made the assumption that $l_j \ne 0$. After implicitly defining the translational and angular velocities of the load with the kinematic equations $\dot{v}_r = x_r$ and $\dot{q}_r = \omega_r \times q_r$, and rearranging terms, we obtain the desired dynamical control system. \hfill$\square$

\subsection{Proof Lemma 1}

\textit{Proof:} Observe that $V_x$ defined as $V_{x} := \frac12 ||e_{v_L}||^2 + c_x e_{x_L}^T e_{v_L}+ \frac12 k_{x_{r}} ||e_{x_L}||^2$ can be bounded from above and below as $\frac12 z_{x}^T \underbar{P}_{x_L} z_{x} \le \ V_{x} \le \frac12 z_{x}^T \bar{P}_{x_L} z_{x}$ where $z_x = \begin{bmatrix} ||e_{x_L}|| & ||e_{v_L}|| \end{bmatrix}^T.$ Further note that both $\underbar{P}_x$ and $\bar{P}_x$ are positive-definite provided that $c_x < \sqrt{k_{x_L}}$.

Similarly, for $j = 1, \cdots, 4$, we define $V_{q_j} := \frac12 ||e_{\omega_j}||^2 + c_q e_{q_j}^T e_{\omega_j} + k_q \Psi_{q_j}$, which is bounded as $\frac12 z_{q_j}^T \underbar{P}_{q_j} z_{q_j} \le V_{q_j} \le \frac12 z_{q_j}^T \bar{P}_{q_j} z_{q_j},$ where $z_{q_j} = \begin{bmatrix} ||e_{q_j}|| & ||e_{\omega_j}|| \end{bmatrix}^T$. As before, $\underbar{P}_{q_j}$ and $\bar{P}_{q_j}$ are positive-definite  when $c_q < \sqrt{k_q}.$

Finally, for $V_{R_L} := \frac12 e_{\Omega_L}^T J_L \Omega_L + k_{R_L} \Psi_{R_L} + c_R e_{R_L}^T J_L e_{\Omega_L}$, we have $\frac12 z_{R_L}^T \underbar{P}_{R_L} z_{R_L} \le V_{R_L} \le \frac12 z_{R_L}^T \bar{P}_{R_L} z_{R_L},$ where $z_{R_L} = \begin{bmatrix} ||e_{R_L}|| & ||e_{\Omega_L}|| \end{bmatrix}^T$, and $\underbar{P}_{R_L}$ and $\bar{P}_{R_L}$ are positive-definite  when $c_{R_L} < \frac{\sqrt{2k_{R_L}\lambda_{\min}(J_L)}}{\lambda_{\max}(J_L)}.$

Observing that $V = V_x + \sum_{j = 1, \cdots, 4} V_{q_j} + V_{R_L}$, we then have that our Lyapunov candidate is bounded as $\frac12 z^T \underbar{P} z \le V \le \frac12 z^T \bar{P} z,$ where
\begin{align*}
    z = \big[ ||e_{x_L}|| \quad ||e_{v_L}|| \quad ||e_{R_L}|| \quad &||e_{\Omega_L}|| \quad ||e_{q_1}|| \quad \cdots \quad ||e_{q_4}|| \\ &||e_{\omega_1}|| \quad \cdots \quad ||e_{\omega_4}|| \big]^T,
\end{align*} and $\bar{P} = \bar{P}_{x_L} \oplus \bar{P}_{R_L} \oplus \bar{P}_{q_1} \oplus \cdots \oplus \bar{P}_{q_4}$, 
    $\underbar{P} =\underbar{P}_{x_L} \oplus  \underbar{P}_{R_L} \oplus \underbar{P}_{q_1} \oplus \cdots \oplus \underbar{P}_{q_4}$.

Next, note that by the invariance of circular shifts of the scalar triple product and the fact that $q_j^T e_{q_j} = 0$, we have:
\begin{align*}
\frac{d}{dt}\Psi_{q_j} &= -\tilde{q}^T_j \dot{q}_j - q_j^T \dot{\tilde{q}}_j = -\tilde{q}_j^T (\omega_j \times q_j) - q_j^T (\tilde{\omega}_j \times \tilde{q}_j)\\ 
&= \omega_j^T (\tilde{q}_j \times q_j) - \tilde{\omega}_j^T (\tilde{q}_j \times q_j) = (\omega_j - \tilde{\omega}_j)^T e_{q_j}\\ &= (\omega_j + (q_j^T \tilde{\omega}) q_j - \tilde{\omega}_j)^T e_{q_j} = e_{\omega_j}^T e_{q_j}.
\end{align*}
Additionally, from the vector triple product, we see
\begin{align*}
\dot{e}_{q_j} &= (\dot{\tilde{q}}_j \times q_j) + (\tilde{q}_j \times \dot{q}_j) = (\tilde{\omega}_j \times \tilde{q}_j) \times q_j - (\omega_j \times q_j) \times \tilde{q}_j \\
&= (\omega_j + \tilde{\omega}_j) \times e_{q_j} + (\tilde{q}_j^T q_j)e_{\omega_j} -  (\tilde{q}_j^T q_j)  (q_j^T \tilde{\omega}_j)q_j \\
&= e_{\omega_j} \times e_{q_j} + (\tilde{q}_j^T q_j)e_{\omega_j} + 2\tilde{\omega}_j \times e_{q_j} -  (\tilde{q}_j^T q_j)  (q_j^T \tilde{\omega}_j)q_j,
\end{align*}
so that
\begin{align*}
\dot{e}_{q_j}^T e_{\omega_j} &= (\tilde{q}_j^T q_j)||e_{\omega_j}||^2 + 2(\tilde{\omega}_j \times e_{q_j})^T e_{\omega_j}\\ &\le ||e_{\omega_j}||^2 + C_{q_j} ||e_{q_j}|| ||e_{\omega_j}||,
\end{align*} \normalsize where $C_{q_j} \le 2\sup ||\tilde{\omega}_j||$ is a non-negative constant. Therefore, the time derivative of the proposed Lyapunov function is bounded as
 \begin{align*}
\dot{V}  \le& -(k_{v_L} - c_x) ||e_{v_L}||^2 + c_x k_{v_L} ||e_{v_L}|| ||e_{x_L}|| - c_xk_{x_L}||e_{x_L}||^2 \\&+ (||e_{v_L}|| + c_x ||e_{x_L}||) ||Y_x|| \\
&(k_{\Omega_L} - 2c_R \lambda_{\max}(J))||e_{\Omega_L}||^2 + c_R (k_{\Omega_L} + B) ||e_{R_L}|| ||e_{\Omega_L}|| \\
&-c_R k_{R_L} ||e_{R_L}||^2 + (e_{\Omega_L} + c_R e_{R_L})||Y_R||\\
&- c_q k_q ||e_{q_j}||^2 -(k_{\omega} - c_q) ||e_{\omega_j}||^2 + c_q (k_{\omega} + C_{q_j} ) ||e_{q_j}|| ||e_{\omega_j}||,
\end{align*} where $Y_x$ and $Y_R$ satisfy the inequalities

\begin{align*}
    (&||e_{v_L}|| + c_x ||e_{x_L}||) ||Y_x|| \le \\
    &\sum_{j=1}^4 \alpha_j \beta(c_x k_{x_L} ||e_{x_L}||^2 + c_x k_{v_L} ||e_{x_L}|| ||e_{v_L}|| + k_{v_L} ||e_{v_L}||^2 ) \\
    & \big[c_x B||e_{x_L}|| + (\beta k_{x_L} e_{x_{\max}} + B)||e_{v_L}||\big]||e_{q_j}|| \\
    &+\alpha_j \gamma(c_x ||e_{x_L}|| + ||e_{v_L}||)(k_{R_L}||e_{R_L}|| + k_{\Omega_L}||e_{\Omega_L}||)
\end{align*}
and
\begin{align*}
    (&||e_{\Omega_L}|| + c_R ||e_{R_L}||) ||Y_R|| \le \\
    &\sum_{j=1}^4 \alpha_j \sigma_j(c_R k_{R_L} ||e_{R_L}||^2 + c_R k_{\Omega_L} ||e_{R_L}|| ||e_{\Omega_L}|| + k_{\Omega_L} ||e_{\Omega_L}||^2 ) \\
    & \big[c_R B||e_{R_L}|| + (\alpha_L \sigma_j k_{R_L} + B)||e_{\Omega_L}||\big]||e_{q_j}|| \\
    &+\alpha_j \delta_j \gamma(c_R ||e_{R_L}|| + ||e_{\Omega_L}||)(k_{x_L}||e_{x_L}|| + k_{v_L}||e_{v_L}||)
\end{align*}
Applying these inequality directly to the above bound for $\dot{V}$, we find that $\dot{V} \le -\sum_{j=1}^4 z_j^T W_jz_j$, where $z_j = \begin{bmatrix} ||e_{x_L}|| & ||e_{v_L}|| & ||e_{R_L}|| & ||e_{\Omega_L}|| & ||e_{q_j}|| & ||e_{\omega_j}|| \end{bmatrix}$. Hence, if each $W_j$ is positive-definite, it follows that the origin is an exponentially stable equilibrium point. $\hfill\square$

\printbibliography

@article{PalCruIRAM12,
	Author = {I. Palunko and P. Cruz and R. Fierro},
	Journal = {IEEE Robotics and Automation Magazine},
	Number = {3},
	Pages = {69--79},
	Title = {Agile Load Transportation},
	Volume = {19},
	Year = {2012}}

@article{rey_dynamics_1999,
	title = {Dynamics and control of a class of underactuated mechanical systems},
	volume = {44},
	issn = {1558-2523},
	pages = {1663--1671},
	number = {9},
	journaltitle = {{IEEE} Transactions on Automatic Control},
	author = {Reyhanoglu, M. and van der Schaft, A. and Mcclamroch, N.H. and Kolmanovsky, I.},
	year = {1999},
}

@inproceedings{bloch2017variational,
  title={Variational obstacle avoidance problem on Riemannian manifolds},
  author={Bloch, Anthony and Camarinha, Margarida and Colombo, Leonardo},
  booktitle={2017 IEEE 56th Annual Conference on Decision and Control (CDC)},
  pages={145--150},
  year={2017},
  organization={IEEE}
}

@article{bloch2021dynamic,
  title={Dynamic interpolation for obstacle avoidance on Riemannian manifolds},
  author={Bloch, Anthony and Camarinha, Margarida and Colombo, Leonardo J},
  journal={International Journal of Control},
  volume={94},
  number={3},
  pages={588--600},
  year={2021},
  publisher={Taylor \& Francis}
}

@article{colombo2015variational,
  title={Variational integrators for mechanical control systems with symmetries},
  author={Colombo, Leonardo and Jim{\'e}nez, Fernando and de Diego, David Mart{\'\i}n},
  journal={Journal of Computational Dynamics},
  volume={2},
  number={2},
  pages={193},
  year={2015},
  publisher={American Institute of Mathematical Sciences}
}

@article{colombo2016geometric,
  title={Geometric integrators for higher-order variational systems and their application to optimal control},
  author={Colombo, Leonardo and Ferraro, Sebasti{\'a}n and de Diego, David Mart{\'\i}n},
  journal={Journal of Nonlinear Science},
  volume={26},
  number={6},
  pages={1615--1650},
  year={2016},
  publisher={Springer}
}

@article{lee2017global,
  title={Global formulations of Lagrangian and Hamiltonian dynamics on manifolds},
  author={Lee, Taeyoung and Leok, Melvin and McClamroch, N Harris},
  publisher={Springer}
}

@article{kobilarov2011discrete,
  title={Discrete geometric optimal control on Lie groups},
  author={Kobilarov, Marin B and Marsden, Jerrold E},
  journal={IEEE Transactions on Robotics},
  volume={27},
  number={4},
  pages={641--655},
  year={2011},
  publisher={IEEE}
}

@inproceedings{kobilarov2014discrete,
  title={Discrete optimal control on lie groups and applications to robotic vehicles},
  author={Kobilarov, Marin},
  booktitle={2014 IEEE international conference on robotics and automation (ICRA)},
  pages={5523--5529},
  year={2014},
  organization={IEEE}
}

@article{izadi2014rigid,
  title={Rigid body attitude estimation based on the Lagrange--d’Alembert principle},
  author={Izadi, Maziar and Sanyal, Amit K},
  journal={Automatica},
  volume={50},
  number={10},
  pages={2570--2577},
  year={2014},
  publisher={Elsevier}
}

@article{lee_feedback_2009,
	title = {Feedback linearization vs. adaptive sliding mode control for a quadrotor helicopter},
	volume = {7},
	issn = {1598-6446},
	pages = {419--428},
	number = {3},
	journaltitle = {International Journal of Control, Automation and Systems},
	author = {Lee, Daewon and Jin Kim, H. and Sastry, Shankar},
	date = {2009-06-01},
	langid = {english},
}

@inproceedings{raffo_backsteppingnonlinear_2008,
	title = {Backstepping/nonlinear {$H_\infty$} control for path tracking of a quadrotor unmanned aerial vehicle},
	pages = {3356--3361},
	booktitle = {Proceedings of the American Control Conference},
	author = {Raffo, Guilherme V. and Ortega, Manuel G. and Rubio, Francisco R.},
	year = {2008},
}

@inproceedings{frazzoli_trajectory_2000,
	title = {Trajectory tracking control design for autonomous helicopters using a backstepping algorithm},
	volume = {6},
	pages = {4102--4107},
	booktitle = {Proceedings of the American Control Conference},
	author = {Frazzoli, E. and Dahleh, M.A. and Feron, E.},
	year = {2000},
}

@inproceedings{sreenath_geometric_2013,
	title = {Geometric control and differential flatness of a quadrotor {UAV} with a cable-suspended load},
	eventtitle = {52nd {IEEE} Conference on Decision and Control},
	pages = {2269--2274},
	booktitle = {52nd {IEEE} Conference on Decision and Control},
	author = {Sreenath, Koushil and Lee, Taeyoung and Kumar, Vijay},
	date = {2013-12},
	note = {{ISSN}: 0191-2216},
}

@inproceedings{lee_geometric_2010,
	title = {Geometric tracking control of a quadrotor {UAV} on {SE}(3)},
	eventtitle = {49th {IEEE} Conference on Decision and Control ({CDC})},
	pages = {5420--5425},
	booktitle = {49th {IEEE} Conference on Decision and Control ({CDC})},
	author = {Lee, Taeyoung and Leok, Melvin and {McClamroch}, N. Harris},
	date = {2010-12},
	note = {{ISSN}: 0191-2216},
}

@article{sreenath1dynamics,
  title={Dynamics, Control and Planning for Cooperative Manipulation of Payloads Suspended by Cables from Multiple Quadrotor Robots},
  author={Sreenath, Koushil and Kumar, Vijay},
  journal={rn},
  volume={1},
  number={r2},
  pages={r3}
}

@inproceedings{lee_geometric_2013,
	title = {Geometric control of cooperating multiple quadrotor {UAVs} with a suspended payload},
	eventtitle = {52nd {IEEE} Conference on Decision and Control},
	pages = {5510--5515},
	booktitle = {52nd {IEEE} Conference on Decision and Control},
	author = {Lee, Taeyoung and Sreenath, Koushil and Kumar, Vijay},
	date = {2013-12},
	note = {{ISSN}: 0191-2216},
}

@inproceedings{lee_geometric_2014,
	title = {Geometric control of multiple quadrotor {UAVs} transporting a cable-suspended rigid body},
	eventtitle = {53rd {IEEE} Conference on Decision and Control},
	pages = {6155--6160},
	booktitle = {53rd {IEEE} Conference on Decision and Control},
	author = {Lee, Taeyoung},
	date = {2014-12},
	note = {{ISSN}: 0191-2216},
}

@inproceedings{wu_geometric_2014,
	title = {Geometric control of multiple quadrotors transporting a rigid-body load},
	eventtitle = {53rd {IEEE} Conference on Decision and Control},
	pages = {6141--6148},
	booktitle = {53rd {IEEE} Conference on Decision and Control},
	author = {Wu, Guofan and Sreenath, Koushil},
	date = {2014-12},
	note = {{ISSN}: 0191-2216},
}

@article{lee_geometric_2018,
	title = {Geometric Control of Quadrotor {UAVs} Transporting a Cable-Suspended Rigid Body},
	volume = {26},
	issn = {1558-0865},
	pages = {255--264},
	number = {1},
	journaltitle = {{IEEE} Transactions on Control Systems Technology},
	author = {Lee, Taeyoung},
	date = {2018-01},
	note = {Conference Name: {IEEE} Transactions on Control Systems Technology},
}

@article{khalil,
	title = {Nonlinear Systems},
	journaltitle = { Prentice-Hall, Inc., Third Edition},
	author = {H. K. Khalil},
	date = {2002},}

@article{khalil2,
	title = {Nonlinear Control},
	journaltitle = {Pearson Education Limited},
	author = {H. K. Khalil},
	date = {2015},}

@article{gas,
	title = {Dynamic collaboration without communication: Vision-based cable-suspended load transport with two quadrotors},
	journaltitle = {IEEE International Conference on Robotics and Automation},
	author = {M. Gassner, T. Cieslewski, D. Scaramuzza},
	pages = {5196--5202},
	date = {2017}}

@article{HSS,
	title = {Geometric mechanics and symmetry},
	journaltitle = {Oxford University Press},
	author = {Holm, D. D and Schmah, T and Stoica, C},
	date = {2009},}

@article{abloch,
	title = {Nonholonomic mechanics and control},
	journaltitle = { Springer-Verlag, Second Edition},
	author = {Bloch A. M},
	date = {2015},}

@article{MR,
	title = {Introduction to Mechanics and Symmetry: A Basic Exposition of Classical Mechanical Systems},
	journaltitle = { Springer-Verlag},
	author = {Marsden, J. E and  Ratiu, T. S},
	date = {1999},}

@article{Leebook,
	title = {Global formulations of Lagrangian and Hamiltonian dynamics on manifolds},
	journaltitle = {Springer 13:31},
	author = {Lee, T and Leok, M and McClamroch. N. H},
	date = {2017},}

@article{GF,
	title = {Calculus of variations},
	journaltitle = {Revised English edition translated and edited by Richard A. Silverman. Prentice-Hall, Inc., Englewood Cliffs, N.J.,},
	author = {Gelfand, L and Fomin, S},
	date = {1963},}

@online{amazonair,
  author = {Amazon.com Inc.},
  title = {Amazonair},
  year = 2021,
  url = {https://www.amazon.com/Amazon-Prime-Air/b?ie=UTF8&node=8037720011},
  urldate = {2021-10-20}
}

@online{googlewing,
  author = {X Development LLC},
  title = {X-{W}ing},
  year = 2021,
  url = {https://x.company/projects/wing/},
  urldate = {2021-10-20}
}

@article{bayen2006adjoint,
  title={Adjoint-based control of a new {E}ulerian network model of air traffic flow},
  author={Bayen, Alexandre M and Raffard, Robin L and Tomlin, Claire J},
  journal={IEEE Transactions on Control Systems Technology},
  volume={14},
  number={5},
  pages={804--818},
  year={2006},
  publisher={IEEE}
}

@article{chung2018survey,
  title={A survey on aerial swarm robotics},
  author={Chung, Soon-Jo and Paranjape, Aditya Avinash and Dames, Philip and Shen, Shaojie and Kumar, Vijay},
  journal={IEEE Transactions on Robotics},
  volume={34},
  number={4},
  pages={837--855},
  year={2018},
  publisher={IEEE}
}

@incollection{ribeiro2021multi,
  title={Multi-robot Systems for Precision Agriculture},
  author={Ribeiro, Angela and Conesa-Mu{\~n}oz, Jesus},
  booktitle={Innovation in Agricultural Robotics for Precision Agriculture},
  pages={151--175},
  year={2021},
  publisher={Springer}
}

@article{god,
title = {Dynamics and control of quadrotor uavs transporting a rigid body connected via flexible cables},
journaltitle = {in 2015 American Control Conference,},
author = {F. A. Goodarzi and T. Lee},
pages = {4677--4682},
date = {2015}}

@article{jo,
title = {Position control of a flexible cable gantry
crane: theory and experiment},
journaltitle = {in American Control Conference,
Proceedings of the 1995,},
author = {S. Joshi and C. D. Rahn},
pages = {2820--2824},
date = {1995}}

@article{jacobquads,
	title = {Geometric Control of two Quadrotors Carrying a Rigid Rod with Elastic Cables},
	journaltitle = {arXiv preprint arXiv:2104.06155,},
	author = {Goodman, J. R and Colombo, L. J},
	date = {2021},}

@article{SPT,
	title = {Introduction to perturbation methods},
	journaltitle = {Springer Science $\&$ Business Media},
	author = {M. Holmes},
	date = {2012},}

@inproceedings{sreenath_geometric_2017,
	title = {Dynamics and control of a quadrotor with a payload suspended
through an elastic cable},
	eventtitle = {52nd {IEEE} Conference on Decision and Control},
	pages = {3906--3913},
	booktitle = {2017 American Control Conference},
	author = {Kotaru, P and Wu, G and Sreenath, K},
	date = {2017-09},
}

\end{document}